\documentclass[aos]{imsart}

\RequirePackage{amsthm,amsmath,amsfonts,amssymb}
\RequirePackage[numbers]{natbib}
\RequirePackage[colorlinks,citecolor=blue,urlcolor=blue]{hyperref}
\RequirePackage{graphicx, bm}
\usepackage{todonotes}

\newcommand{\norm}[1]{\left\lVert#1\right\rVert}

\newcommand{\Db}{\mathbf{D}}
\newcommand{\Eb}{\mathbf{E}}

\newcommand{\Hb}{\mathbf{H}}

\newcommand{\Pb}{\mathbf{P}}

\newcommand{\bA}{\bm{A}}

\newcommand{\bC}{\bm{C}}

\newcommand{\bP}{\bm{P}}

\newcommand{\bS}{\bm{S}}

\newcommand{\bU}{\bm{U}}

\newcommand{\bX}{\bm{X}}

\newcommand{\bbeta}{\bm{\beta}}
\newcommand{\bgamma}{\bm{\gamma}}

\newcommand{\cL}{\mathcal{L}}

\newcommand{\cN}{\mathcal{N}}

\newcommand{\bb}{\bm{b}}

\newcommand{\bp}{\bm{p}}

\newcommand{\bv}{\bm{v}}

\newcommand{\bx}{\bm{x}}

\newcommand{\RR}{\mathbb{R}}

\newcommand*{\zero}{{\bm 0}}

\newcommand{\logit}{{\rm logit}}

\DeclareMathOperator{\Var}{{\rm Var}}

\DeclareMathOperator*{\Cov}{\rm Cov}

\newcommand{\argmin}{\mathop{\mathrm{argmin}}}

\newcommand{\tr}{\mathop{\mathrm{tr}}}

\newcommand\independent{\protect\mathpalette{\protect\independenT}{\perp}}
\def\independenT#1#2{\mathrel{\rlap{$#1#2$}\mkern2mu{#1#2}}}

\newcommand{\indep}{\rotatebox[origin=c]{90}{$\models$}}

\startlocaldefs
\theoremstyle{plain}

\newtheorem{theorem}{Theorem}[section]
\newtheorem{lemma}[theorem]{Lemma}
\newtheorem{corollary}{Corollary}[section]
\theoremstyle{remark}

\newtheorem{remark}{Remark}[section]

\endlocaldefs

\begin{document}

\begin{frontmatter}
\title{BELIEF in Dependence: Leveraging Atomic Linearity\\ in Data Bits for Rethinking Generalized Linear Models}
\runtitle{BELIEF in Dependence}

\begin{aug}
\author[A]{\fnms{Benjamin}~\snm{Brown}\ead[label=e1]{brownb1@live.unc.edu}},
\author[A]{\fnms{Kai}~\snm{Zhang}\ead[label=e2]{zhangk@email.unc.edu}}
\and
\author[B]{\fnms{Xiao-Li}~\snm{Meng}\ead[label=e3]{meng@stat.harvard.edu}}
\address[A]{Department of Statistics and Operations Research, University of North Carolina at Chapel Hill\printead[presep={,\ }]{e1,e2}}

\address[B]{Department of Statistics, Harvard University \printead[presep={,\ }]{e3}}
\end{aug}

\begin{abstract}
Two linearly uncorrelated binary variables must be also independent because non-linear dependence cannot manifest with only two possible states. This inherent linearity is the atom of dependency constituting any complex form of relationship. Inspired by this observation, we develop a framework called binary expansion linear effect (BELIEF) for understanding arbitrary relationships with a binary outcome. Models from the BELIEF framework are easily interpretable because they describe the association of binary variables in the language of linear models, yielding convenient theoretical insight and striking Gaussian parallels. With BELIEF, one may study generalized linear models (GLM) through transparent linear models, providing insight into how the choice of link affects modeling. For example, setting a GLM interaction coefficient to zero does not necessarily lead to the kind of no-interaction model assumption as understood under their linear model counterparts. Furthermore, for a binary response, maximum likelihood estimation for GLMs paradoxically fails under complete separation, when the data are most discriminative, whereas BELIEF estimation automatically reveals the perfect predictor in the data that is responsible for complete separation. We explore these phenomena and provide related theoretical results. We also provide preliminary empirical demonstration of some theoretical results.
\end{abstract}



\end{frontmatter}

\section{Nonparametric Modeling Through Data Bits}\label{sec: intro}

\subsection{Taking Advantage of an Inherent Linearity}

There are two kinds of classical scientists: those who believe their models, and those who model their belief. As such, misspecification is the fundamental gremlin of statistical modeling: incorrect models tempt practitioners with mathematical elegance while quietly belying reality. Statistics and data science have long sought modeling strategies free of unduly restrictive assumptions, spanning from traditional settings \citep[e.g.,][]{mccullagh2019generalized,htf2009} to more recent efforts \citep[e.g.,][]{lei2018distribution,buja2019,barber2020distribution,gupta2020distribution,foygel2021limits,li2020nonparametric}. 

This paper reports an effort to derive a general modeling theory using the framework of binary expansion statistics (BEStat) \citep{zhang2019bet,zhang2021beauty}, which leads to multi-resolution linear models for a binary outcome under a fully nonparametric setting. Without any assumption on their joint distribution, random variables can be effectively decomposed into data bits. These data bits can be regarded as the \textit{atoms} of information from both statistical and computer science perspectives. By constructing models and formulating inference directly from the bits, this framework provides additional theoretical insight on the binary world and suggests a new approach for analyzing generalized linear models (GLMs), as well as an alternative modeling strategy.

To understand complex forms of dependency, we begin by studying the simplest form of dependency---the dependency between two binary variables. Consider two $\pm 1$-valued random variables $A$ and $B$. Trivially, because $\Eb[B|A]$ can only take two states depending on the value of $A$, it follows that
\begin{equation}\label{eq: belief1}
   \Eb[B|A]= \beta_0 + \beta_1 A, 
\end{equation}
which is intrinsically a linear model with slopes $\beta_0, \beta_1\in\mathbb{R}$. In the special case that $A$ and $B$ are Rademacher, $\Eb[B] = 0$ necessitates that $\beta_0 = 0$, while $\Var(A) = 1$ implies that $\beta_1 = \Eb[AB] = \Cov(A, B)/\Var(A) $. More generally, suppose $\Pb(A = 1) = \theta$, $\Pb(B = 1\mid A = 1) = \phi$, and $\Pb(B = 1\mid A = -1) = \psi$. Here, one can show $\bbeta_0 = \phi + \psi - 1$, while 
\begin{equation}
\bbeta_1 = \phi - \psi = \frac{4\theta(1-\theta)(\phi - \psi)}{4\theta(1-\theta)} = \frac{\Cov(A, B)}{\Var(A)}.
\end{equation}
In any case, since the conditional distribution $\Pb(B|A)$ is determined by its mean given in \eqref{eq: belief1}, $A$ and $B$ are independent if and only if $\bbeta_1 = 0$. The linearity in this atomic case inspires us to think about the possibility of modeling any form of association through binary variables.

Linearity is inherent in the binary nature of $A$ and $B$, because any non-linearity requires more than two states to reveal. {We emphasize that the key observation of linearity is completely general.} It is useful to recognize that this atomic linearity {carries} to arbitrarily many binary predictors to construct a saturated linear model of the conditional probability. Moreover, the slopes, which we term BELIEF coefficients {in the context of a binary probability model}, are unique when the second moment matrix of the binary predictors is positive definite and can be estimated through the least squares algorithm, as we show in Theorem~\ref{thm:sigma}. {Whenever this uniqueness does not hold, we say that the predictors have a degenerate distribution, and we discuss the situation in depth in Section~\ref{subsec: singular}.}

Besides indicating the strength and direction of the dependence, the slopes from this intrinsic linearity are also useful in specifying dependence structures in the joint distribution. For example, there is a direct correspondence  between the conditional independence and multiplicative subgroups of binary predictors with nonzero slopes, as shown in Theorem~\ref{thm: subgroup}. This connection provides a basis for regularization in estimation and prediction problems, where screening of variables \citep{fanlv2008,he2013quantile,zhang2017spherical} or interactions \citep{fan2015innovated, thanei2018xyz} are well understood. Moreover, the boundedness of binary variables and their slopes (Theorem~\ref{thm:sigma}) facilitates the applications of machine learning and high-dimensional statistical methods \citep{buhlmann2011}. 

Unlike the simple linear model in \eqref{eq: belief1}, a GLM connects slopes with the response by way of a (typically nonlinear) link function. These slopes are estimated through maximum likelihood estimation (MLE), and the specific choice of link function in the likelihood dictates how a practitioner interprets model coefficients. For instance, in the famed logistic model for a binary response, a slope describes the effect of the corresponding predictor on the log-odds. In contrast, interpretation of BELIEF coefficients occurs directly on the level of cell probabilities. Well-known in econometrics, a related model called the linear probability model (LPM) achieves similar interpretability by expressing the conditional response probability directly as a linear function of the predictors. However, this approach does not consider the binary expansion of explanatory variables. Indeed, for sufficiently extreme predictor values, the response probability falls outside $[0, 1]$, which is recognized as a substantial drawback compared to GLMs \citep{angrist2009mostly, wooldridge2010econometric}. Via binary expansion into data bits, BELIEF achieves linearity while also guaranteeing valid response probabilities. The recent \cite{vansteelandt2022assumption} investigates assumption-lean inference for main effects and interactions in GLMs, describing how these statistical properties may still be captured even in the case of model misspecification. We study the connections between GLMs and BELIEF in Section~\ref{subsec: generalresultsforglms}.

In addition, the linearity illustrated in \eqref{eq: belief1} is evocative of classical Gaussian linear models. The Gaussian analogy turns out to be a recurring theme in the study of binary random variables, as the binary world shares many familiar Gaussian properties, while differing in unintuitive ways. Table~\ref{tab: compare} summarizes the comparison of properties between binary and Gaussian variables that we will establish in this article.

\begin{table}
\caption{A comparison of distributional and inferential properties for binary and Gaussian linear models.}
\begin{center}
\begin{tabular}{c c|l }
 Binary & Gaussian & Property \\ 
 \hline
  Y & Y &Independence is equivalent to uncorrelation (Theorem \ref{thm:sigma})\\
  Y & Y & Conditional expectation is a linear equation with slopes $\bbeta$ (Theorem \ref{thm:sigma})\\
  Y & Y & Slopes that are zero relate to conditional independence (Theorem \ref{thm: subgroup})\\
 Y & Y & Least squares $\hat{\bbeta}$ is the MLE (asymp. normal + efficient) (Theorem \ref{thm:asymnorm})\\
 Y & Y & $\hat{\bbeta}$ is unbiased whenever it exists (Theorem \ref{thm:whenMPunbiased})\\
  N & Y &  Existence of residuals that are independent of $\hat{\bbeta}$ in general (supplementary material \citep{beliefsupp})\\
  N & Y & KL-divergence is free of marginal information, in a sense described in \cite{beliefsupp}\\
  Y & N & Slopes are within a compact convex set (Theorem~\ref{thm:sigma}).
\end{tabular}
\end{center}
\label{tab: compare}
\end{table}

Going beyond binary variables, we can use binary expansion to approximate uniform variables to an arbitrary accuracy, and hence to approximate any (continuous) variable via the probability integral transformation. This fact is summarized in the following lemma in \cite{zhang2021beauty}.

\begin{lemma}\label{lem: binary:approx} 
Let $\bU=(U_1,U_2,\cdots, U_p)^\top$ be a random vector supported within $[-1,1]^p$. There exists a sequence of binary random variables $\{A_{j,d}\}$, $j=1,2,\cdots,p$, $d=1,2,\cdots,D$, which take only values $-1$ and $1$, such that $
\max_{1\le j\le p}\{ |U_j  - U_{j,D}|\} \to 0$ almost surely as $D \rightarrow \infty$, where $U_{j,D} = \sum_{d=1}^D \left(A_{j,d}\right)/2^d.$
\end{lemma}

By Lemma~\ref{lem: binary:approx}, for any random variable $U$, the first $D$ data bits $\{A_d\},d=1,\ldots D,$ form a filtration  to approximate the distribution of $U,$ and  $\sigma_D=\sigma(A_1,\ldots,A_D)$ is the $\sigma$-field summarizing all information up to depth $D$ in the binary expansion. Hence $D$ is a resolution level, as in the multi-resolution framework of \cite{li2021multi}.  Because of the aforementioned inherent linearity, when $U$ is used as a predictor for a binary response, there is an intrinsic equation expressing the conditional expectation of the response as a linear function of the binary variables in $\sigma_D$. We are thus able to approximate the dependency between $U$ and the binary response by extracting the hidden linearity through the binary expansion approach. 

By combining the linear dependency of binary variables and the binary expansion approximation of the distribution, there is a general distribution-free modeling strategy built upon atomic linearity in data bits, as we show in Section~\ref{subsec: doublelimit}. We thus refer to this modeling framework as the \textit{binary expansion linear effect} (BELIEF).

\subsection{Revisiting GLMs with BELIEF} \label{subsec: illustrative}

For the better part of a century, GLMs have been a prevalent tool for modeling  binary outcome, as summarized in  \cite{mccullagh2019generalized}. Popular methods such as logistic and probit regressions work well when class probabilities are monotone in the predictors but struggle otherwise. The log-linear model is another useful GLM for contingency tables, where the linearity is an assumption over log cell probabilities. 

To showcase the BELIEF framework as applied to GLMs, we begin with an illustrative example. Let $A_1, A_2$, and $B$ be binary random variables taking values $\pm 1$. A GLM model linking $B$ to $A_1, A_2$ is given by 
\begin{equation}\label{eq: glmtoy}
\Pb(B = 1|A_1, A_2) = g^{-1}(\gamma_0 + \gamma_1A_1 + \gamma_2A_2 + \gamma_{12}A_1A_2),
\end{equation}
where $g$ is a chosen link function. The familiar expression in \eqref{eq: glmtoy} represents a saturated GLM, in the sense that all possible predictor interactions are included, as in for example \cite{bahadur1961representation,Wahba1994,dai2013multivariate}.

Alternatively, without any link function or assumption, one can show 
\begin{equation}\label{eq: belieftoy}
\Pb(B = 1|A_1, A_2) = \beta_0 + \beta_1A_1 + \beta_2A_2 + \beta_{12}A_1A_2
\end{equation}
for some slopes $\bbeta = (\beta_0, \beta_1, \beta_2, \beta_{12})^\top$. Here, $\bbeta$ is constrained to $\norm{\Hb_4\bbeta}_{\infty} \leq 1$, where $\Hb_4$ is the $4\times 4$ Hadamard matrix according to Sylvester's construction. The equation \eqref{eq: belieftoy} expresses the conditional distribution of $B$ as a multilinear function of the predictors directly. We refer to $\bbeta$ as the vector of BELIEF coefficients, in contrast to the GLM coefficients $\bgamma$ for link function $g$.

At heart, our main insights for GLMs rely on the comparison between \eqref{eq: glmtoy} and \eqref{eq: belieftoy}. Globally, neither \eqref{eq: glmtoy} nor \eqref{eq: belieftoy} restricts the joint distribution of $(A_1, A_2, B)$. As long as $(A_1, A_2, B)$ takes the maximum of $2^3 = 8$ possible values with positive probability, then for any joint distribution on $(A_1, A_2, B)$, there exists exactly one $\bgamma$ satisfying \eqref{eq: glmtoy} and exactly one $\bbeta$ satisfying \eqref{eq: belieftoy}, forming a bijection $\bbeta \leftrightarrow \bgamma$.

However, the story changes for connecting individual components of $\bgamma$ with that of $\bbeta$. For example, the null hypothesis $H_0: \gamma_{12} = 0$ is not equivalent to $H_0: \beta_{12}=0$, that is, there is no interaction effect between $A_1$ and $A_2$ governing the conditional distribution of $B$. Furthermore, no choice of link function resolves this issue---the statement $(\gamma_{12} = 0 \iff \beta_{12} = 0)$ implies that $g$ itself is linear, as discussed in Section \ref{subsec: generalresultsforglms}. Indeed, the unintuitive relationship between link functions and interaction terms is known to be problematic in practice \citep{berry2010testing, rainey2016compression}.

Such differences between $\bbeta$ and $\bgamma$ can go only so far, as both vectors must respect probabilistic properties that are invariant to the model representation. For example, $\beta_1 = \beta_2 = 0$ if and only if $\gamma_1 = \gamma_2 = 0$, because both conditions are equivalent to $B \independent{(A_1, A_2)} | A_1A_2$. In Section~\ref{sec: belief}, we will show in Theorem~\ref{thm: subgroup} that these essential similarities can be captured by a group structure on the predictors. 

We would like to remark that binary variables are often coded as $0/1$ in traditional GLM settings. Here, we adopt the Rademacher $-1/1$ for mathematical convenience. With this arithmetically symmetric coding, multiplication of two binary variables corresponds to a group operation on the set of binary variables and their interactions, which we relate to statistical modeling properties. Moreover, this coding does not sacrifice any generality in the breadth of binary models that may be represented. For general statistical interpretations, the choice of coding will be content dependent; see \cite{cox1972analysis} and \cite{mccullagh2000invariance}.

Aside from its theoretical value as a tool for model comparison, BELIEF suggests a unified modeling strategy in its own right. Since the conditional expectation of a binary outcome determines its distribution, BELIEF enjoys many theoretically optimal properties by directly estimating this crucial quantity with least squares. It also avoids an inherent defect with GLM, which leads to non-identifiability of $\bbeta$ when the prediction problem is the easiest, that is, when there is a prefect predictor \citep{albert1984existence}, as studied in Section~\ref{subsec: separation}. However, there is no free lunch. BELIEF, as a modeling strategy, suffers significantly more from the curse of dimensionality compared to GLM, because the BELIEF model parameters grow exponentially in the number of predictors. This is especially the case for using binary expansions to approximate continuous predictors as in Lemma~\ref{lem: binary:approx}. Nevertheless, the dimensionality of the BELIEF framework faithfully reflects the inherent challenges of distribution-free modeling for arbitrary relationships---where the complexity of the problem naturally grows exponentially. Despite these difficulties, BELIEF remains valuable as a general approach for understanding the optimal trade-offs among key factors such as sample size, number of predictors, strength of relationships, and relationship complexity. In Section \ref{subsec: scalability}, we propose techniques to facilitate greater practical applicability. In the supplementary materials \citep{beliefsupp}, we provide an impossibility result articulating the trade-off between dimensionality and sensitivity to interaction. We will study these topics more intensively in future work.

\subsection{Additional Relatives of BELIEF}\label{subsec: related}
The filtration approach in BELIEF can be viewed as a ``probably approximately correct'' learning strategy articulated by \cite{valiant1984theory} in computer science. It has connections to several important modeling methods. As we were reminded by a reviewer, because our approach comes with guarantees based on mathematical proof, it might be better termed ``pro\textbf{v}ably approximately correct'' instead.

Machine learning algorithms, such as tree-based methods, are an important approach for classification with binary outcomes \citep{htf2009}. The discretization achieved by the binary expansion can also be viewed as a tree-based information decomposition strategy, but with more principled and theoretically tractable constructions and justifications. For instance, the orthogonality among the bits eliminates the redundancy in information bases, such as those obtained by greedy--partitioning tree constructions. Furthermore, the linearity of BELIEF provides clear interpretation of the effect from each bit. This linearity has been studied in the analysis of Boolean functions $f:\{-1, 1\}^p \rightarrow \{-1, 1\}$. For example, \cite{o2014analysis} addresses the representation of Boolean functions, albeit leading to different goals than we have here, given our own emphasis on statistical modeling. In addition, \cite{vanderweele2008empirical} discusses the interpretation of binary linear model parameters in the context of sufficient cause interactions.

The binary expansion of marginal variables is also closely related to wavelet analysis \citep{antoniadis2012wavelets} and the tensor product models \citep{Wahba1994,lin2000}; the former decomposing random variables instead of their functions, as done in the latter frameworks. The advantage of direct stochastic decomposition/representation is that it enhances interpretability by explicating functional relationships, which can better facilitate model developments and computational implementations. Another related linear model of dependence between binary variables is the Ising model \citep{ravikumar2010high,chen2016fused}. 

Recent developments on distribution-free modeling of dependency include \cite{lei2018distribution}, \cite{gupta2020distribution} and \cite{foygel2021limits}. An important recent result is the impossibility of constructing a confidence interval for the label probability to have uniform coverage for all distributions of the covariates \citep{barber2020distribution}. This type of non-uniformity result is also common in testing problems, for example in \cite{zhang2019bet}. The binary expansion framework avoids such non-uniformity issues by restricting the inference up to certain meaningful resolution levels. We explain the details of the inference in Section~\ref{sec: estimate}.

Furthermore, BELIEF offers a different discretization strategy from estimating proportions of the binary response in each cell after a naive discretization of covariates. As an analogy, the connection between the two approaches is similar to the estimation in the physical domain and the frequency domain in Fourier or wavelet analysis \citep{antoniadis2012wavelets}. By studying the dependency over the binary variables as a basis, interpretations on the global forms of dependency can be learned, as described in the testing problems \citep{zhang2019bet,zhang2021beauty}. Moreover, the principled discretization incorporated with resolution provides meaningful information on the approximation of the dependence. In particular, the sparsity of slopes is closely connected to conditional independence. Aside from studies of (in)dependence, binary transformations of arbitrary random variables are also sufficient to demonstrate ``association'' in the specific sense of \cite{esary1967association}, which is a stronger condition than nonnegative correlation, but neither implies nor forbids dependence.

The rest of the article is organized as follows: Section~\ref{sec: belief} introduces BELIEF and discusses characterizations of correlation and independence for binary models. Section~\ref{sec: estimate} considers estimation properties of BELIEF. Section~\ref{sec: rethinking} provides insights from BELIEF regarding GLMs and other models of dependence. Section~\ref{sec: improvementAndPractical} evaluates extensions and practical concerns. Section~\ref{sec: empirical} illustrates the method with numerical studies. Section~\ref{sec: discussion} discusses directions for further explorations.

\section{Basic Theory of the BELIEF Framework}\label{sec: belief}

\subsection{The Inherent Linearity in Binary Variables} \label{sec: inherent}
Now, we introduce essential notation for BELIEF. Let $\bm{A} = (A_1, A_2, \dots, A_p)^\top\in \{-1, 1\}^p$ be a binary random vector,  and define the corresponding ``power vector" of dimension $2^p$
by taking the Kronecker product:
\begin{align*}
\bm{A}^{\otimes} := \begin{pmatrix} 1 \\ A_1 \end{pmatrix} \otimes \begin{pmatrix} 1 \\ A_2 \end{pmatrix} \otimes \cdots \otimes \begin{pmatrix} 1 \\ A_p \end{pmatrix}.
\end{align*} 
That is, $\bm{A}^{\otimes}$ collects all possible combination of components of $\bm{A}$, taken as a product. As shown in \cite{zhang2019bet}, the $2^p$ possible states of $\bm{A}^{\otimes}$ form the Hadamard matrix $\Hb_{2^p}$ according to Sylvester's construction.

The simplicity of the binary world implies natural properties of the $\sigma$-field generated by a binary random vector $\bA.$ One of the most important such inherent properties is the linearity of measurable functions. The linearity is summarized in the following two lemmas, which also demonstrate the critical roles of $\bm{A}^{\otimes}$. 
\begin{lemma}\label{lem: pd}
	For any binary vector $\bA=(A_1,\ldots,A_p)^\top$, the following properties are equivalent for its probability space.
	\begin{enumerate}
		\item $\Eb[\bA^{\otimes}\bA^{\otimes^\top}]$ is positive definite.
		\item The vector $\bA$ takes $2^p$ distinct values all with positive probabilities.
		\item The $\sigma$-field $\sigma(A_1,\ldots,A_p)$ contains $2^{2^p}$ distinct events all with positive probabilities, except for the empty set.
	\end{enumerate}
\end{lemma}
\begin{lemma}\label{lem: linearity}
Let $\bA$ be a binary vector of length $p$. For any function $f:\{-1, 1\}^p \rightarrow \mathbb{R}$, there exists $\bbeta\in\mathbb{R}^{2^p}$ such that $f(\bA) = \bbeta^\top\bA^\otimes$.
\end{lemma}
Lemma~\ref{lem: pd} states the importance of the second moment matrix of binary variables. Lemma~\ref{lem: linearity} further asserts that any function over binary variables is linear. This important feature substantially facilitates the understanding of complicated relationships and simplifies the analysis.

Building on these basic properties, we reveal the beauty of binary random variables manifested as similarities with the Gaussian world. Table \ref{tab: compare} summaries similarities and differences between the binary and Gaussian worlds that are examined in this article. In particular, dependence of $\bm{A}$ on another binary random vector $\bm{B}$ is characterized by coordinate-wise dependence of $\bm{A}^{\otimes}$ and $\bm{B}^{\otimes}$. Moreover, for a single binary random variable $B$ distributed jointly with $\bm{A}$, its conditional expectation given $\bm{A}$ is a linear function of $\bm{A}^{\otimes}$ because $\bm{A}^{\otimes}$ provides a full basis to express all $2^p$ possible values of $\Eb[B|\bm{A}]$ as a function of $\bm{A}$. This linear relationship between $B$ and $\bm{A}^{\otimes}$ is the essence and technical meaning of
BELIEF. 
\setcounter{theorem}{0}
\begin{theorem}\label{thm:sigma}
With the notation above, the following properties hold:
\begin{enumerate}\label{th:betabound}
\item The binary vectors $\bm{A} = (A_1, \dots, A_{p_1})^\top$ and $\bm{B} = (B_1,\dots , B_{p_2})^\top$ are independent if and only if ${\bm{A}}^{\otimes}$ and ${\bm{B}}^{\otimes}$ are uncorrelated for all the coordinate pairs $({\bm{A}}^{\otimes}_i,{\bm{B}}^{\otimes}_j)$.

\item For a binary variable $B \in \{-1,1\}$ and a $p$-dimensional binary vector $\bm{A}\in \{-1, 1\}^p$, we have 
\begin{equation}\label{eq: belief}
    \Eb[B|\bm{A}] = \bm{\beta}^\top\bm{A}^{\otimes}
\end{equation} for some constant vector $\bm{\beta}\in\mathbb{R}^{2^p}$.

\item $\bbeta = \Eb^{-1}[\bm{A}^{\otimes}{\bm{A}^{\otimes}}^\top]\Eb[\bm{A}^{\otimes} B]$ whenever the inverse of $\Eb[\bA^{\otimes}\bA^{\otimes^\top}]$ exists. In this case, $\bbeta$ satisfies
(a) $|\bbeta^\top\bA^{\otimes}| \leq 1$, (b) $\norm{\Hb_{2^p}\bbeta}_\infty \leq 1$, and (c) $\norm{\bbeta}_2 \leq 1$, where $\Hb_{2^p}$ is the Hadamard matrix of size $2^p$ according to Sylvester's construction.
\end{enumerate}
\end{theorem}

Parallel to the famous Gaussian property, Part 1 of Theorem \ref{thm:sigma} reduces the assessment of independence of binary variables $\bm{A}$ and $\bm{B}$ to checking pairwise correlations of binary variables in $\sigma(\bm{A})$ and $\sigma(\bm{B})$. Parts 2 and 3 of Theorem \ref{thm:sigma}, which comprise the heart of BELIEF, are a natural consequence of Lemma~\ref{lem: linearity}, as $\Eb[B|\bA]$ is a measurable function of $A_1,\ldots,A_p$. While the mathematical foundations of BELIEF are straightforward, the statistical implications are comparatively subtle. Classical textbooks teach us to study the relationship between $B$ and $\bA$ through a logistic regression, a probit regression, or other forms of GLM. However, without any model assumptions, Theorem~\ref{thm:sigma} asserts that we are able to achieve optimal inference of dependence simply through least squares of data bits. Here the optimality is in terms of the mean squared error (MSE), as seen by $\Eb[B|\bA]=\argmin_f{\Eb[(B-f(\bA))^2]}.$ 

Here the role of $\bbeta$ in a linear combination of predictors is reminiscent of Gaussian conditional expectation, and rightly so---the closed-form expression for $\bbeta$ is strikingly familiar.
Breaking slightly from the Gaussian parallels, it is worth noting that the boundedness of binary random variables leads to corresponding bounds on $\bbeta$, stated in part 3 of Theorem \ref{thm:sigma}. Unlike the normal family, the family of binary vectors do not permit location shift since it is not a location family, leading to a more restrictive parameter space. Nonetheless, in the study of dependence, the sign and magnitude of each slope indicates the direction and strength of relationship between the binary interaction and the response, as in classical linear models. Moreover, the hypercubic parameter space for BELIEF slopes with ``corners'' hints at natural methods of regularization with least squares optimization, which we discuss in Section~\ref{subsec: scalability}. Vertices of this hypercube correspond to models in which $\bA$ determines the value of $B$, and in the context of estimation, these corners characterize the phenomenon of perfect separation (Section~\ref{subsec: separation}).

Furthermore, the fact that $\Pb(B=r|A)$ is fully characterized by $\Eb[B| A]$ via
\begin{align}\label{eq:cond}
\Pb(B = r| A) = \frac{1+r\Eb[B| A]}{2},  \quad {\rm for}\quad r\in\{1, -1\}
\end{align}
gives a special interpretation to coefficients in $\bbeta$ that are zero. Note that the conditional probability $\Pb(B = r | A)$ plays crucial roles in many important problems such as the Bayes rule in classification \citep{htf2009} and the propensity score in causal inference \citep{rosenbaum1983central}. Statistically optimal estimation of BELIEF and its efficient least squares algorithm can be helpful in advancing these areas.

Pushing the Gaussian connection farther, the following result is analogous to the well-known result in the Gaussian world, where preservation of coefficients comes with orthogonality.  

\begin{theorem}
\label{thm: slopesAsDepthChanges}
Suppose $A_p$ is a Rademacher random variable independent of $\bA_{-p} = (A_1, \dots, A_{p-1})$. Then for any coordinate $C$ of $\bA_{-p}^\otimes$, the coefficient of $C$ in $\Eb[B|\bA]$ is equal to the coefficient of $C$ in $\Eb[B|\bA_{-p}]$.
\end{theorem}

Theorem \ref{thm: slopesAsDepthChanges} comes directly from an application of the tower property of conditional expectations. More broadly, the technique of writing $\Eb[B|\bA_{-p}] = \Eb[\Eb[B|\bA]|\bA_{-p}]$ can be used to derive general relationships between coefficients as the number of predictors varies.

Notably, BELIEF slopes $\bbeta$ are uniquely defined when $\Eb[\bm{A}^{\otimes}{\bm{A}^{\otimes}}^\top]>0$. According to Lemma \ref{lem: pd}, this condition fails exactly when $\bA$ has empty cells with probability one. For example, let $\bA = (A_1, A_2, A_3)^\top$, wherein $A_3 = A_1A_2$, and observe that $\Pb(\bA = (-1, -1, -1)^\top) = 0$. Importantly, however, the BELIEF slopes in $\Eb[B|\bA]$ may still be uniquely defined when the joint distribution of $(B, \bA^\top)^\top$ has cells with zero probability---uniqueness of the slopes depends on the marginal distribution of the predictors. We discuss the issue of degeneracy in more detail in Section~\ref{subsec: singular}.

\subsection{Algebraic Considerations for Dependence with BELIEF}\label{subsec: condind}

In the Gaussian world, parameters of linearity such as slopes and correlations have important connections to conditional dependence structure. Sparsity in linear models thus plays an important role in modeling. For instance, it is well known that conditional independence induces decomposition of likelihood through a graphical model \citep{lauritzen1996graphical}. As BELIEF makes explicit, dependence in the binary world is also linear, and analogous properties govern the conditional independence of binary variables. However, a proper treatment of these properties must respect algebraic considerations that do not apply in the Gaussian case. For example, given univariate Gaussian variables $X$ and $Y$, the product $XY$ is not generally Gaussian. In contrast, the class of $\pm 1$ binary random variables is closed under multiplication; this principle appears in our construction of the object $\bA^\otimes$, which is a binary vector containing all distinct products of components of $\bA$. 

In probabilistic terms, the algebraic structure of binary variables characterizes sub-$\sigma$-fields generated by their interactions. As a simple example, if we know $A_1A_2$ and $A_2A_3$, then we can deduce the value of $A_1A_3 = A_1A_2\cdot A_2A_3$, but not $A_1$ or $A_3$. We find basic group theory to be a natural way of expressing this structure. Given a set $S$ of binary variables, let $\langle S \rangle$ be the group generated by $S$ under multiplication. Further, let $G_{\bA} = \langle \{A_1, A_2, \dots, A_p\}\rangle$. That is, the components of the vector $\bA^\otimes$ are exactly the elements of the set $G_{\bA}$. Finally, to ensure that group elements are well-defined, we stipulate that the vector $\bA$ satisfy the equivalent conditions of Lemma \ref{lem: pd}, meaning that $\bA$ takes $2^p$ distinct values with positive probability.

\begin{theorem}\label{thm: subsigma}
Let $A_1,\ldots,A_p$ be binary variables with $\Eb[\bA^{\otimes}\bA^{\otimes^\top}]>0$. Every $\sigma$-field generated by a subset $S \subseteq G_{\bA}$ is generated by a unique set of binary variables $H\subseteq G_{\bA}$ that forms a multiplicative group, and $|H| = 2^k$ for some $k$. In particular, $H = \langle S \rangle$, and $H$ is isomorphic to the $k$-factor direct product $\mathbb{Z}_2^k.$
\end{theorem}

Theorem \ref{thm: subsigma} gives basic intuition on the origin and structure of these subgroups. The $k$-factor direct product $\mathbb{Z}_2^k$ refers to the direct product of the additive group $\mathbb{Z}_2$ with itself $k$ times, that is $\mathbb{Z}_2 \times \cdots \times \mathbb{Z}_2$. Crucially, there may be many different subsets of $G_{\bA}$ that generate $\sigma (S)$, but $\langle S \rangle$ is the only one that forms a group under multiplication. While Theorem \ref{thm: subsigma} concerns multiplicative generation of binary variables, the next result comments on the group properties of binary variables from the perspective of data reduction.

\begin{theorem}\label{thm: groupdatareduction}
Let $A_1,\ldots,A_p$ be binary variables with $\Eb[\bA^{\otimes}\bA^{\otimes^\top}]>0$. Let $H$ be a subgroup of $G_{\bA}$, and let $k$ be such that $|H|=2^k$. Then there exists $S\subseteq G_{\bA}$ with $|S| = k$, $\langle S \rangle = H$, and the entries of $S$ jointly take $2^k$ distinct values with positive probability. Further, for any $T\subset G_{\bA}$ with $|T| < k$, we have $\sigma (T) \neq \sigma (H)$.
\end{theorem}

Intuitively, a group of $2^k$ binary variables can be represented by a minimum of $k$ bits. To demonstrate Theorems \ref{thm: subsigma} and \ref{thm: groupdatareduction}, we provide a concrete example. Define $\bA = (A_1, A_2, A_3)^\top$, and let $S = \{ A_1, A_1A_2, A_2A_3\}$ and $T = \{ A_1A_2, A_2A_3\}$. Taking products, one can compute that $\langle S \rangle = G_{\bA}$, while $\langle T \rangle = \{1, A_1A_2, A_1A_3, A_2A_3\}$. By Theorem \ref{thm: subsigma}, we know that $\sigma (S) = \sigma (G_{\bA} ) \neq \sigma (\langle T \rangle )= \sigma (T)$, as well as $G_{\bA} \cong \mathbb{Z}_2 \times \mathbb{Z}_2 \times \mathbb{Z}_2$ and $\langle T \rangle \cong \mathbb{Z}_2 \times \mathbb{Z}_2$. Moreover, by Theorem \ref{thm: groupdatareduction}, we know that $S$ and $T$ are of minimal size for generating their respective groups and $\sigma$-fields.

Now, we are ready to state conditional independence results pertaining to BELIEF slopes $\bbeta$.

\begin{theorem}\label{thm: subgroup}
Let $B,A_1,\ldots,A_p$ be binary variables with $\Eb[\bA^{\otimes}\bA^{\otimes^\top}]>0$ and $\Eb[B|\bA]=\bbeta^\top\bA^{\otimes}.$ 
\begin{enumerate}
    \item Let $\cL=\{\Lambda: \beta_\Lambda \neq 0\}$ be the collection of elements in $G_{\bA}$ whose slopes in $\bbeta$ are nonzero, and let $G_\cL=\langle\cL\rangle$ under multiplication of elements. Then 
	\begin{equation}
	B\independent{G_{\bA}} | G_{\cL}.
	\end{equation}
	\item Let $S \subseteq G_{\bA}$. If $B\independent{G_{\bA}}|S$, then the coefficients of $G_{\bA}\setminus \langle S \rangle$ in $\bbeta$ are zero. 
\end{enumerate}
\end{theorem}

Theorem \ref{thm: subgroup} explains how sparsity of BELIEF slopes governs conditional independence of binary variables. In Part 1 of Theorem \ref{thm: subgroup}, nonzero slopes determine a $\sigma$-field that yields conditional independence of $B$ and $\bA$. Conversely, in Part 2, conditional independence forces certain slopes to be zero. 

We would like to remark that Part 2 of Theorem \ref{thm: subgroup} is not true of $G_{\bA}\setminus S$. For instance, consider $\bA = (A_1, A_2, A_3)^\top$, and $S = \{A_1A_2, A_2A_3\}$. Given $B$ such that $\Eb[B|\bA] = 0.25A_1A_2 + 0.25A_2A_3 + 0.25A_1A_3$, then $B$ is independent of $\bA$ given $S$, but there are nonzero slopes outside $S$. Hence, Gaussian intuition is not able to identify these slopes correctly, indicating the advantage of group notation in the binary setting.

In the context of Theorem \ref{thm: subsigma}, it is reasonable to ask why we bother to condition on $G_{\cL}$ in Part 1 of Theorem \ref{thm: subgroup}, rather than just $\cL$. After all, the distinction is merely formal, for we know $\sigma (G_{\cL}) = \sigma (\cL)$. However, by writing $G_{\cL}$, we intend to emphasize that conditioning on $\textit{any}$ generator of $G_{\cL}$ trivially gives the same conditional independence. Potential generators may exclude some variables whose slopes are nonzero. Conversely, removing a zero slope variable from a generator may break the conditional independence property. This flexibility in representation signifies a substantial departure from typical Gaussian conditional independence.

To elucidate Theorem \ref{thm: subgroup}, we list two concrete examples of conditional independence: 

(1) Consider BELIEF of a binary response $B$ on binary covariates $A_1$ and $A_2$. If $\Eb[B|A_1,A_2]=\alpha+\beta A_1A_2,$ then $B\independent{(A_1,A_2)} | A_1A_2.$ 

(2) Consider BELIEF of a binary response $B$ on binary covariates $A_1$, $A_2,$ and $A_3$. If $\Eb[B|A_1,A_2,A_3]=\alpha+\beta_1 A_1A_2+ \beta_2 A_2A_3+\beta_3 A_1A_3,$ then $B\independent{(A_1,A_2,A_3)} \mid  |C|$, where $C=\frac{A_1}{2}+\frac{A_2}{4}+\frac{A_3}{8}.$ This is because $(A_1A_2, A_2A_3, A_1A_3)$ generates the same $\sigma$-field as $|C|$, as both are solely determined by whether each pair of bits has the same sign. In describing the form of dependence, this means the conditional distribution of $B$ given $C$ is symmetric w.r.t. the origin.

Finally, while we express Theorem \ref{thm: subgroup} using BELIEF coefficients $\bbeta$, analogous properties hold for GLM slopes $\bgamma$ in the context of saturated binary predictors.
\setcounter{theorem}{0}
\begin{corollary}\label{cor: glmsubgroup}
Let $B,A_1,\ldots,A_p$ be binary variables with $\Eb[\bA^{\otimes}\bA^{\otimes^\top}]>0$ and $\Eb[B|\bA]=g^{-1}(\bgamma^\top\bA^{\otimes}).$ 
\begin{enumerate}
    \item Let $\cL=\{\Lambda: \beta_\Lambda \neq 0\}$ be the collection of elements in $G_{\bA}$ whose slopes in $\bgamma$ are nonzero, and let $G_\cL=\langle\cL\rangle$ under multiplication of elements. Then 
	\begin{equation}
	B\independent{G_{\bA}} | G_{\cL}.
	\end{equation}
	\item Let $S \subseteq G_{\bA}$. If $B\independent{G_{\bA}}|S$, then the coefficients of $G_{\bA}\setminus \langle S \rangle$ in $\bgamma$ are zero. 
\end{enumerate}
\end{corollary}
\setcounter{theorem}{6}
Any binary linear model, generalized or not, must respect the connection between group structure and conditional independence. Hence, whenever GLM slopes $\bgamma$ and BELIEF slopes $\bbeta$ describe the same model, $\bgamma$ cannot differ from $\bbeta$ in ways that contradict the group-wise properties of this section. Nonetheless, these properties still permit important structural differences between the slope vectors. In Section \ref{subsec: generalresultsforglms}, we explain how each slope in $\bgamma$ encodes different information than the corresponding slope in $\bbeta$. For example, $\bgamma_i = 0$ does not in general imply $\bbeta_i = 0$, nor the converse.

In summary, the general results above provide an easy way to check for nontrivial conditional independence by noting whether the binary variables with nonzero BELIEF slopes are within a proper subgroup of $G_{\bA}$. {Another perspective on these statements pertains to sufficiency: together with $B$, the set of interactions with nonzero coefficients is sufficient for $\bbeta$.} These results also provide insights for developments of regularization methods. The boundedness of binary variables also make many results in high-dimensional statistics and machine learning applicable \citep{buhlmann2011}.

A special case of these theorems is summarized below in Theorem \ref{thm:condind}, which draws as close a comparison as possible to well-known Gaussian results. As usual, we take $\bA$ to be nondegenerate.

\setcounter{theorem}{5}
\begin{theorem}
\label{thm:condind}
For any $1 \leq i \leq p$, let $J_i$ be the set of indices of interactions in $\bA^{\otimes}$ that involve the binary variable $A_i$. Define the sub-vector $\bA^{\otimes}_{0}$ by removing the leading one from $\bA^{\otimes}$, and let $\bC = (B, \bA^{\otimes \top}_{0})^\top$. If the vector of BELIEF slopes $\bbeta$ is uniquely defined and $\norm{\Hb_{2^p}\bbeta}_\infty < 1$, then the following statements are equivalent:
\begin{enumerate}
\item $\bbeta_{j} = 0$ for all $j\in J_i$.
\item $B$ and $A_i$ are conditionally independent given $(\bA_{-i})^{\otimes}$, where $\bA_{-i}$ is a sub-vector of $\bA$ by removing $A_i$.
\item $[\text{Cov}(\bC)^{-1}]_{j, 1} = 0$ for all $j \in J_i$.
\end{enumerate}
\end{theorem}

The properties above closely mirror conventional knowledge about Gaussian regression slopes and precision matrices, while accounting for the additional complications of binary variables. We further explore theoretical properties of binary variables in the supplementary materials \citep{beliefsupp}.

\section{Estimation with BELIEF}\label{sec: estimate}
\subsection{Least Squares Estimation of BELIEF Slopes}
\label{subsec: leastsquares}
Suppose we observe an i.i.d. sample $\{B_i,  \bA_i\}_{i = 1}^n$, where the predictor $\bA_i$ is a $p$-dimensional binary vector, and the response $B_i$ is binary as well. Let $\mathcal{B}=(B_1, \ldots,  B_n)^\top$, and define a $n\times 2^p$ matrix $\mathcal{A}=(\bA^{\otimes}_1, \ldots, \bA^{\otimes}_n)^\top$, using the notation from Section \ref{sec: belief}. A natural estimate of $\bbeta$ given in Theorem~\ref{th:betabound} is the least squares estimator (LSE) $
\hat{\bbeta} = (\mathcal{A}^\top\mathcal{A})^{-1}\mathcal{A}^\top\mathcal{B},
$ which satisfies the same bounds as the estimand $\bbeta$, as one would expect.

\begin{theorem}
\label{thm:bounds} Whenever the LSE $\hat{\bbeta}$ exists, it satisfies 
(1) $\norm{\Hb_{2^p}\hat{\bbeta}}_\infty \leq 1$, (2) $\norm{\mathcal{A}\hat{\bbeta}}_\infty \leq 1$, and (3) $\norm{\hat{\bbeta}}_2 \leq 1$.
\end{theorem}

This seemingly trivial result is at the heart of BELIEF estimation, because it shows that we can deliver legitimate probability estimates that are linear in
$\hat\beta$, that is, being always between $[0, 1]$,  without ever imposing GLM-type of link function for the same purpose. Consequently, BELIEF eliminates the need to make up model assumptions for mathematical conveniences. This is because, in analogy to \eqref{eq:cond}, the corresponding LSE for $\Pb(B = r|\bA^{\otimes} = \bv)$ is 
\begin{align*}
\hat{\Pb}(B = r|\bA^{\otimes} = \bv) = \frac{1 + r\bv^\top\hat{\bbeta}}{2}, \qquad {\rm for}\ r=-1, 1,
\end{align*}
which never goes outside of the unit interval because $|\bv^\top\hat\bbeta| \leq 1$, by Theorem~\ref{thm:bounds}.

As a caveat, Theorem~\ref{thm:bounds} does not say that performing least squares on arbitrary binary predictors must always yield valid cell probability estimates. In fact, for some binary data sets, naively removing columns from $\mathcal{A}$ results in $\norm{\Hb_{2^p}\hat{\bbeta}}_\infty > 1$. In order to reduce the dimension of the model while preserving the result of Theorem~\ref{thm:bounds}, one must ensure that the remaining binary predictors form a multiplicative group as described in Theorem~\ref{thm: subsigma}. This preservation is due to Theorem~\ref{thm: groupdatareduction}: a subgroup $H\leq G_{\bA}$ of order $2^k$ is generated by some $S\subseteq G_{\bA}$ with $|S| = k$. Letting $\bA'$ denote a $k$-length vector whose components are the $k$ distinct elements of $S$, we see that Theorem~\ref{thm:bounds} applies to the LSE for the reduced BELIEF model $\Eb[B|\bA'] = \bbeta^\top (\bA')^\otimes$. 

To proceed to additional results about $\hat{\bbeta}$, we reiterate the conditions under which the parameter $\bbeta$ is well-defined. It follows from Theorem \ref{thm:sigma} that $\bbeta$ is uniquely defined when $\Eb[\bA^{\otimes}{\bA^{\otimes}}^\top]>0$, or equivalently when $\Pb(\bA^{\otimes} = \bv) > 0$ for all rows $\bv$ of $\Hb_{2^p}$.

\begin{theorem}
\label{thm:asymnorm}
Assume $|\Eb[B|\bA]| < 1$, and $\bbeta$ is uniquely defined. Then the least squares estimator $\hat{\bbeta}$ is also the unique MLE whenever $\hat{\bbeta}$ exists, where the likelihood is with respect to the joint distribution of the $(p+1)\times 1$ binary vector $(B, \bA^\top)^\top$. As $n\rightarrow\infty$, we have $\hat{\bbeta} \xrightarrow{a.s.} \bbeta$ and 
$$\sqrt{n}(\hat{\bbeta} - \bbeta) \xrightarrow{d} N\left(0,\  \frac{1}{2^{2p}}\Hb_{2^p} \Db \Hb_{2^p}\right),$$ 
where 
$$
\Db = \text{diag}\left(\frac{1 - (\bbeta^\top\bv_i)^2}{\Pb(\bA^{\otimes} = \bv_i)}\right)_{i = 1}^{2^p}, $$
and $\bv_i$ is the $i$th row of $\Hb_{2^p}$, $1 \leq i \leq 2^p$. 
In particular, the variance of each entry of $\hat{\beta}$ is the same,
$\tr(D)/2^{2p}$. 
\end{theorem}

Theorem \ref{thm:asymnorm} tells us that LSE is also the MLE under the joint model, which might come as a surprise since the LSE does not use any information about the likelihood function of $\bbeta$. The reason that this is true is given by Theorem \ref{thm:ancillary}, which reveals that $\bA$ is ancillary to $\bbeta$. This is again in parallel to the Gaussian world, or more broadly to the world of regression with continuous variables. Yet, it is not a forgone conclusion, for in the discrete world, marginal distributions often have more impact on the joint distribution than their continuous counterparts; see for example, the discussion of this matter in the context of the $2\times 2$ model \citep[e.g.,][]{chernoff2004information}. 

Relevant to prediction, Theorem \ref{thm:asymnorm} gives the asymptotic variance of the estimate $\hat{\Eb}[B|\bA = \bv_i]$. To achieve the classical binomial confidence interval for each cell, we need only write the asymptotic variance in terms of conditional probabilities. Intuitively, this relationship is because the estimation of each coordinate of $\bbeta$ uses information from all $2^p$ cells, whereas only a single cell is pertinent to each estimated cell expectation in the absence of additional assumptions on the regression function. This is the price paid for extreme modeling flexibility. In Section \ref{subsec: scalability}, we explore potential remedies to BELIEF's high-dimensional issues.

Returning to ancillarity, define $\bp_{\bA^{\otimes}}\in\mathbb{R}^{2^p}$ such that the $i$th component of $\bp_{\bA^{\otimes}}$ is $\Pb(\bA^{\otimes} = \bv)$, where $\bv$ is the $i$th row of $\Hb_{2^p}$. One can think of $\bp_{\bA^{\otimes}}$ as containing all information about the marginal distribution of $\bm{A}$.

\begin{theorem}
\label{thm:ancillary}
If $\bbeta$ is uniquely defined, the joint distribution of $(\bA, B)$ is characterized by $(\bp_{\bA^{\otimes}}, \bbeta)$, and $\bA$ is ancillary for $\bbeta$. In addition, $\bp_{\bA^{\otimes}}$ and $\bbeta$ are orthogonal parameters.
\end{theorem}

Theorem \ref{thm:ancillary} states that we can treat the relationship between $\bA$ and $B$ as an entirely separate matter from the distribution of $\bA$ itself. The orthogonality of $\bp_{\bA^{\otimes}}$ and $\bbeta$ means that, in terms of asymptotic efficiency, it makes no difference whether or not we know the marginal parameters $\bp_{\bA^{\otimes}}$. 

Theorem~\ref{thm:asymnorm} permits us to construct asymptotic confidence intervals for components of $\hat{\bbeta}$ or predictions $\bv^\top\hat{\bbeta}$. We can also perform inference on the BELIEF regression function: for each dyadic cell, we can give a confidence interval for the probability that $B = 1$.

Inference aside, in Theorem \ref{thm:asymnorm}, we mention the fact that all the slopes in $\hat{\beta}$ have the same asymptotic variance. These equal variances facilitate regularization---there is no need to adjust the penalty according to the variance of each slope. This fact, as the result of equal aggregations of variance through an orthogonal transformation with $\Hb_{2^p}$, represents an advantage over a naive discretization of covariates.

As a final point, we state a convergence property relevant to dimension $p\rightarrow \infty$, using results given by \cite{chernozhukov2022high}. In contrast to Theorem \ref{thm:asymnorm}, we assume nonstochastic covariates, and we impose structural assumptions to prevent degeneracy and satisfy moment conditions. Let $p$ represent the number of binary predictors, such that the parameter space has dimension $2^p$. Let $\{B_{i}\}_{i = 1}^n$ be independent binary random variables such that $\Eb[B_{i}] = \bbeta^\top\bv_{i}$, where $\bv_{i}$ is some row of $\Hb_{2^{p}}$ and $\bbeta$ is some vector in $\mathbb{R}^{2^{p}}$. Define $n_{\bv} = \#\{i: \bv = \bv_{i}\}$. That is, $n_{\bv}$ is the count of predictor observations in the $\bv$ cell. Let $\mathcal{A}$ be the $n\times 2^{p}$ matrix whose $i$th row is $\bv_{i}$, and let $\mathcal{B}$ be the vector of length $n$ whose $i$th entry is $B_{i}$. Now, assume that

(1) there exists some $\epsilon > 0$ such that $\text{Var}(B_{i}) \geq \epsilon$ for all $1 \leq i \leq n$, and

(2) there exist some $0 < K_1 < K_2$ such that $K_1/2^{p} \leq n_{\bv}/n \leq K_2/2^{p}$ for all rows $\bv$ of $\Hb_{2^{p}}$.

\noindent Given that (2) holds, it follows that $\hat{\bbeta}_n = (\mathcal{A}^\top\mathcal{A})^{-1}\mathcal{A}^\top\mathcal{B}$ is well-defined, and we may denote $\bS_n := \sqrt{n}(\hat{\bbeta}_n - \bbeta)$ and $\Sigma := \text{Var}(\bS_n)$. Finally, adapting notation from \cite{chernozhukov2022high}, define
\begin{align*}
\mathcal{R} = \left\{\prod_{j = 1}^{2^p} [a_j, b_j]: -\infty \leq a_j \leq b_j \leq \infty, 1 \leq j \leq 2^p\right\}.
\end{align*}

\begin{theorem}\label{thm: highdim}
If (1) and (2) are satisfied, then
\begin{align*}
\sup_{R\in\mathcal{R}} \left|\Pb(\bS_n \in R) - \Pb(N(0, \Sigma)\in R)\right| \leq  C\frac{\log^{5/4}(2^{p}n)}{n^{1/4}},
\end{align*}
where $C$ is a constant depending on $\epsilon$, $K_1$, and $K_2$.
\end{theorem}

Assumption (1) ensures that the response variables do not have variance tending to zero, and assumption (2) asserts that the proportion of observations in each cell are not imbalanced in the limit. Rather than giving a particular limiting distribution, Theorem \ref{thm: highdim} uniformly bounds the error in normal approximation for the rescaled BELIEF slope estimate $\bS_n$, making the dependence on $p$ explicit. We hope to continue our study of BELIEF in growing dimensions in future work.

\subsection{Dealing with Singularity}\label{subsec: singular}

Though we have much convenient theory regarding $\hat{\bbeta}$, this estimator has drawbacks. The estimate $\hat{\bbeta}$ exists only when $\mathcal{A}^\top\mathcal{A}$ is invertible, which is true exactly when we have at least one observation for each possible value of $\bA^{\otimes}$. When working with the binary expansion of continuous variables, the situation may require estimation at a higher resolution than what least squares can accommodate.

One solution to this invertibility issue is to employ the Moore-Penrose pseudoinverse; that is, we use $\hat{\bbeta}_{MP} := (\mathcal{A}^\top\mathcal{A})^+\mathcal{A}^\top\mathcal{B}$. 
Clearly, $\hat{\bbeta}_{MP} = \hat{\bbeta}$ whenever $\hat{\bbeta}$ exists. When the estimator $\hat{\bbeta}$ does not exist, $\hat{\bbeta}_{MP}$ is agnostic in the following sense.

\begin{theorem}
\label{thm:agnostic}
Suppose there is no observation $\bA^{\otimes}_i$ such that $\bA^{\otimes}_i = \bv$. Then using the pseudoinverse estimator $\hat{\bbeta}_{MP}$, we have
$\hat{\Eb}[B|\bA^{\otimes} = \bv] =\bv^\top\hat{\bbeta}_{MP} = 0,$
or equivalently
\begin{align*}
\hat{\Pb}(B = r|\bA^{\otimes} = \bv) & = \frac{1 + r\bv^\top\hat{\bbeta}_{MP}}{2} = \frac{1}{2}, \qquad {\rm for\ any}\  r\in\{-1, 1\}.
\end{align*}
\end{theorem}
\noindent Ultimately, this theorem follows from the fact that such a $\bv$ must be in the nullspace of $\mathcal{A}$. In general, the estimator $\hat{\bbeta}_{MP}$ gives us a way to extend the usual LSE to a form that exists with probability one.

Owing to the fact that the original LSE $\hat{\bbeta}$ exists with probability less than one in finite samples, we are unable to say whether $\hat{\bbeta}$ is unbiased in the conventional sense. As $\hat{\bbeta}_{MP}$ always exists, we can in fact show \citep{beliefsupp} that $\hat{\bbeta}_{MP}$ is generally biased for $\bbeta$. However, this pseudoinverse estimator comes close to unbiasedness in the following sense.

\begin{theorem}
\label{thm:whenMPunbiased}
Let $S = \{\hat{\bbeta}_{MP} = \hat{\bbeta}\}$. If $\bbeta$ is uniquely defined, then $\Eb[\hat{\bbeta}_{MP}|S] = \bbeta$, and $\hat{\bbeta}_{MP}$ maximizes the likelihood function of $\bbeta$ for the distribution of $(B, \bA^\top)^\top.$
\end{theorem}
\noindent Thus, while $\hat{\bbeta}_{MP}$ may be biased in general,  $\hat{\bbeta}_{MP}$ is unbiased whenever the usual LSE $\hat{\bbeta}$ exists. The maximum likelihood property of Theorem \ref{thm:whenMPunbiased} is trivial whenever $\hat{\bbeta}$ exists, as Theorem \ref{thm:asymnorm} indicates that $\hat{\bbeta} = \hat{\bbeta}_{MP}$ must be the unique MLE. Whenever $\hat{\bbeta}$ fails to exist, there are infinitely many slopes maximizing the likelihood, and owing to a property of the pesudoinverse, $\hat{\bbeta}_{MP}$ has minimal $L^2$ norm among these maximizing slopes.

As an assumption to Theorems \ref{thm:asymnorm} and \ref{thm:whenMPunbiased}, we require that the $p$-vector $\bA$ of binary predictors take $2^p$ distinct values with positive probability, which is the minimal requirement for the true slopes $\bbeta$ to satisfy $\Eb[B|\bA] = \bbeta^\top\bA^\otimes$ uniquely. Theorem \ref{thm:whenMPunbiased} deals with degeneracy at the sample level, wherein empty cells render $\mathcal{A}^\top\mathcal{A}$ singular. Instead, if cells of $\bA$ are empty at the model level, then additional assumptions on $\bbeta$ are required to maintain identifiability, as infinitely many vectors $\bb$ will achieve $\Eb[B|\bA] = \bb^\top\bA^\otimes$. For example, one may consider the theoretical analog to $\hat{\bbeta}_{MP}$, which we denote $\bbeta_{MP} = \{\Eb[\bm{A}^{\otimes}{\bm{A}^{\otimes}}^\top]\}^{+}\Eb[\bm{A}^{\otimes} B]$. This parameter enjoys the property that $\bbeta_{MP}^\top \bv = \hat{\bbeta}_{MP}^\top \bv = 0$ for any row $\bv$ of $\Hb_{2^p}$ satisfying $\Pb(\bA^\otimes = \bv) = 0$. In addition, by imitating the proof of Theorem \ref{thm:asymnorm}, it follows that the estimator $\hat{\bbeta}_{MP}$ for the parameter $\bbeta_{MP}$ also maximizes the likelihood function of $\bbeta_{MP}$ for the distribution of $(B, \bA^\top)^\top$ in the model-degenerate case.

Below, we provide a summary of the possible situations one may encounter when working with binary data, as well as how BELIEF responds in each case. We proceed in order of increasing degeneracy:

\begin{enumerate}
\item \textit{$\Pb(\bA^\otimes = \bv) > 0$, $\#\{i: \bA_i = \bv\} > 0$, and $|\Eb[B|\bA^\otimes = \bv]| < 1$ for all rows $\bv$ of $\Hb_{2^p}$.} This is the easy case. LSE $\hat{\bbeta}$ exists without need for the pesudoinverse, and as the sample size $n$ increases, $\sqrt{n}(\hat{\bbeta} - \bbeta)$ is asymptotically normal with a nonsingular limiting covariance matrix.
\item \textit{$\Pb(\bA^\otimes = \bv) > 0$ and $\#\{i: \bA_i = \bv\} > 0$ for all rows $\bv$ of $\Hb_{2^p}$, but $|\Eb[B|\bA^\otimes = \bv']| = 1$ for some $\bv'$.} Here, we have at least one observation for every possible value of $\bA$, but the value of $B$ is deterministic given $\bA^\otimes = \bv'$. We again have that $\hat{\bbeta}$ exists without need for the pesudoinverse, but this time, $\sqrt{n}(\hat{\bbeta} - \bbeta)$ converges to a degenerate normal distribution. This is because $$\hat{\Eb}[B|\bA^\otimes = \bv'] - \Eb[B|\bA^\otimes = \bv'] = \bv'^{\top}(\hat{\bbeta}-\bbeta) = 0$$ whenever $\hat{\bbeta}$ exists.
\item \textit{$\Pb(\bA^\otimes = \bv) > 0$ for all rows $\bv$ of $\Hb_{2^p}$, but $\#\{i: \bA_i = \bv'\} = 0$ for some $\bv'$.} Even though all theoretical cell probabilities are positive, we were unlucky and failed to observe an outcome in one of the cells. As a result, $\hat{\bbeta}$ does not exist, so we resort to $\hat{\bbeta}_{MP}$, which has the property that $\hat{\bbeta}_{MP}^\top\bv' = \hat{\Eb}[B|\bA^\otimes = \bv'] = 0$. Still, as the sample size increases, we will eventually observe some data point in the $\bv'$ cell. Hence, $\sqrt{n}(\hat{\bbeta} - \bbeta)$ is asymptotically normal with nonsingular limiting covariance.
\item \textit{$\Pb(\bA^\otimes = \bv') = 0$ for some row $\bv'$ of $\Hb_{2^p}$.} This is the worst case---no method has any hope of estimating $\Eb[B|\bA^\otimes = \bv']$. The binary probability model as parameterized by BELIEF slopes is unidentifiable without additional regularity conditions. If we introduce the restriction that $\bbeta^\top\bv' = 1/2$ for all $\bv'$ corresponding to degenerate cells, then $\hat{\bbeta}_{MP}$ is consistent for $\bbeta$, and $\sqrt{n}(\hat{\bbeta}_{MP} - \bbeta)$ converges to a degenerate normal distribution.
\end{enumerate}
In general, the problem of singularity is closely related to collinearity of covariates in classical linear models. Therefore, theory and methods treating collinearity, such as regularization approaches, can be borrowed to alleviate such issues.  

\subsection{Continuous Predictors with Unknown Marginal Distributions}
While BELIEF is formulated with predictors in the form of binary variables, we consider continuous variables using their binary expansion. Ideally, given a continuous random variable $X$ with CDF $F$, we would use the rescaled binary digits of $F(X)\sim U[0, 1]$. However, in practice, we may not know the function $F$. In such a case, it is reasonable to consider the empirical CDF $\hat{F}_{\bX}$, where $\bX$ is our collection of independent realizations of $X$. In this case, one may wonder whether estimated BELIEF slopes $\hat{\bbeta}_{ecdf}$ that depend on empirical CDF (ECDF) transformation are still consistent for the true slopes $\bbeta$. The following result removes this concern. 
\begin{theorem}
\label{thm:unknownmargin}
Under the conditions in Theorem~\ref{thm:asymnorm}, $\hat{\bbeta}_{ecdf} \stackrel{a.s.}{\rightarrow} \bbeta$ as $n\rightarrow \infty$.
\end{theorem}

We have further investigated the asymptotic distribution of the BELIEF coefficients when the empirical CDF transformation is applied to continuous predictors. This analysis reveals some interesting non-standard problems on the asymptotic theory of an empirical process. In the following theorem, we present results for the case where there is one continuous predictor, and BELIEF is fit using the first bit of this predictor, with the empirical CDF transformation applied. In particular, we show that under the null hypothesis of independence, the asymptotic variance of the slope estimate remains the same whether the true CDF or the empirical CDF is used for the marginal transformation.

\begin{theorem}
\label{thm:unknownmargin_p1d1}
For an even number $n$, consider $n$ i.i.d. observations $(X_i,B_i), i=1,\ldots,n$ from $(X,B) \in \RR\times\{-1,1\}$. Denote the population median of $X$ by $M_X$ and the sample median of $X_i$'s by $\hat{M}_X.$ Define $\hat{A}_i=2I(X_i>\hat{M}_X)-1, i=1,\ldots,n.$ Consider the BELIEF estimate $\hat{\bbeta}_{ecdf}$ by replacing $A_i$'s in $\hat{\bbeta}$ as in Section~\ref{sec: estimate} with $\hat{A}_i$'s. Then for the slope $\beta$ and its estimate $\hat{\beta}_{ecdf},$ as $n\rightarrow \infty,$
\begin{equation}
		\sqrt{n}(\hat{\beta}_{ecdf}-\beta) \stackrel{d}{\rightarrow}  \cN(0,\sigma_{ecdf}^2) 
\end{equation}
where
\begin{equation*}
	\sigma_{ecdf}^2=\sigma_{cdf}^2+4(\Pb(B=1|X=M_X)-\Pb(B=1))^2
\end{equation*}
with $\sigma_{cdf}^2$ being the asymptotic variance in Theorem~\ref{thm:asymnorm}. In particular, when $X$ and $B$ are independent, $\sigma_{ecdf}^2=\sigma_{cdf}^2.$
\end{theorem}

Theorem~\ref{thm:unknownmargin_p1d1} provides insight into the validity of inference about the slope estimates with empirical CDF transformed data, which we employ in the numerical studies in Section~\ref{sec: empirical}. Due to the complexity and non-standard nature of the general asymptotic theory, we are only able to present results for the simplest case in this paper. The complete asymptotic theory of the BELIEF slopes with empirically CDF transformed data is undoubtedly very important and very interesting. However, as this theory warrants a separate and in-depth investigation, we plan to address it thoroughly in a future paper.

\section{Rethinking GLMs and Other Models of Dependence}
\label{sec: rethinking}

\subsection{The Meaning of Slopes under a Link Function}
\label{subsec: generalresultsforglms}
Here, we provide results connecting GLMs with the BELIEF framework, building on the concepts introduced in Section \ref{subsec: illustrative}. Broadly, our goal is to understand dependency of random variables via the structure inherited from the binary world.

First, we explain how BELIEF and GLMs are equivalent representations when the predictors are binary.

\begin{theorem} \label{thm: 1-1GLM}
Let $g:(-1, 1)\rightarrow\mathbb{R}$ be a link function, and assume that the binary vector $(B, A_1, \dots, A_p)$ takes $2^{p+1}$ distinct values with positive probability. Then there exist unique vectors $\bgamma$ and $\bbeta$ such that $\Eb[B|\bA] = g^{-1}(\bgamma^\top \bA^\otimes) = \bbeta^\top \bA^\otimes$.
\end{theorem}

In general, given a GLM for $B|\bA$ with slopes $\bgamma$, we can always find an equivalent BELIEF representation with some slopes $\bbeta$, and vice versa. Neither model imposes a restriction on the distribution of $B|\bA$, and the distinction is merely formal.

An exception to the correspondence between $\bbeta$ and $\bgamma$ occurs when $|\Eb[B|\bA^\otimes = \bv]| = 1$ for some row $\bv$ of $\Hb_{2^p}$. The statement of Theorem \ref{thm: 1-1GLM} excludes this case by stipulating that the joint vector $(B, \bA^\top)$ be nondegenerate, but this phenomenon of perfect separation is a common issue for GLMs in practice. For instance, if $g$ is the logistic link and $\Eb[B|\bA^\otimes] = I(\bA^\otimes = \bv)$, then there is no $\bgamma\in\mathbb{R}^{2^p}$ such that $\Eb[B|\bA^\otimes] = g^{-1}(\bgamma^\top \bA^\otimes)$, as the right hand side can never equal one. Alternatively, we can see that $\Eb[B|\bA^\otimes] = \frac{1}{2^p}\bv^\top\bA^\otimes$ holds for BELIEF slopes $\frac{1}{2^p}\bv$.

As Theorem \ref{thm: 1-1GLM} suggests, for models expressible in both $\bbeta$ and $\bgamma$, the bijection $\bbeta \leftrightarrow \bgamma$ carries important information relating the two modeling schemes. The following result explores this correspondence locally via Taylor expansion.

\begin{theorem} \label{thm: taylor}
Let $g:(-1, 1)\rightarrow\mathbb{R}$ be a link function with differentiable inverse $g^{-1}$. Consider the GLM expressed by $\Eb_{\bgamma}[B|\bA] = g^{-1}(\bgamma^\top \bA^\otimes)$. Then $\Eb_{\bgamma + \mathbf{\delta}}[B|\bA] = (\bbeta_{\bgamma}^\top + \mathbf{\delta}^\top\bbeta_{\bgamma}')\bA^\otimes + o(\norm{\mathbf{\delta}}_2)$, where $\bbeta_{\bgamma} = \Hb_{2^p}^{-1}g^{-1}(\Hb_{2^p}\bgamma)$ for $g^{-1}$ applied component-wise. Moreover, $\bbeta_{\bgamma}' = 2^{-p}\Hb_{2^p}\Db\Hb_{2^p}$, where $\Db = \text{diag}((g^{-1})'(\bgamma^\top \bv_i))_{i = 1}^{2^p}$ and $\bv_i$ is the $i$th row of $\Hb_{2^p}$.
\end{theorem}

That is, by Taylor expanding $\Eb_{\bgamma + \mathbf{\delta}}[B|\bA]$, we state how the sensitivity of $\bbeta$ to perturbations in $\bgamma$ depends on the derivative of $g$. This theorem often has a particularly convenient form when expanding at $\bgamma = \mathbf{0}$, as we show in the examples below.

(1) When $g:[-1, 1]\rightarrow\mathbb{R}$ is a logistic link, $\Eb_{\mathbf{0} + \mathbf{\delta}}[B|\bA] \approx \frac{1}{2}\mathbf{\delta}^\top\bA^\otimes$.

(2) When $g:[-1, 1]\rightarrow\mathbb{R}$ is a probit link, $\Eb_{\mathbf{0} + \mathbf{\delta}}[B|\bA] \approx \sqrt{\frac{2}{\pi}}\mathbf{\delta}^\top\bA^\otimes$.

These expansions compare the local scale of logistic and probit parameters to their BELIEF counterparts, which are always on the scale of cell expectations. As another perspective on link functions, we may consider algebraic properties, rather than analytic ones. When do two models agree on the presence of a pairwise interaction?

\begin{theorem} \label{thm: interactionNotEquiv}
Let $g:(-1, 1)\rightarrow \mathbb{R}$ be continuous and nondecreasing with inverse $g^{-1}$. Let $\bA$ be a non-degenerate binary vector of length $p$. Define $\bbeta_{\bgamma}$ to be the $2^p$-length vector satisfying $\bbeta_{\bgamma}^\top \bA^\otimes = g^{-1}(\bgamma^\top \bA^\otimes)$, where $\bbeta_{\bgamma}$ is defined for all $\bgamma$ such that $\bgamma^\top\bA^\otimes$ is in the range of $g$. 

Let $\bgamma_{12}$ and $\bbeta_{\bgamma, 12}$ be the coefficients of $A_1A_2$ from $\bgamma$ and $\bbeta_{\bgamma}$, respectively. Suppose that $\bbeta_{\bgamma}$ satisfies $(\bgamma_{12} = 0 \iff \bbeta_{\bgamma, 12} = 0)$. Then $g$ is a linear function.
\end{theorem}

No meaningful choice of $g$ results in equivalent interactions---the only possible options amount to rescalings of BELIEF. Interestingly, Theorem \ref{thm: interactionNotEquiv} can be extended to a comparison of two different GLMs. Given candidate links $h$ and $k$, then $h \circ k^{-1}$ may take the place of $g$ in Theorem \ref{thm: interactionNotEquiv}.

The principle of Theorem \ref{thm: interactionNotEquiv} extends to GLMs with continuous predictors, as well. Consider the simple model
\begin{equation*}
\Eb[B|U_1, U_2] = g^{-1}(\gamma_0 + \gamma_1U_1 +\gamma_2U_2),
\end{equation*}
where $U_1, U_2 \sim \text{Unif}(-1, 1)$ are independent. In the spirit of Lemma \ref{lem: binary:approx}, we can approximate this model to any depth $d$ using the binary expansion of $U_1$ and $U_2$, where the approximation error $R_d$ converges to zero in $d$. (See Section \ref{subsec: doublelimit} for additional discussion on the error $R_d$ of this approximation.)

Taking the binary expansion, we obtain
\begin{equation*} 
\Eb[B|U_1, U_2] = g^{-1}\left(\gamma_0 + \gamma_1\sum_{i = 1}^d \frac{A_{1, i}}{2^i} +\gamma_2\sum_{i = 1}^d \frac{A_{2, i}}{2^i}\right) + R_d.
\end{equation*}
Letting $\bA = (A_{1, 1}, \dots, A_{1, d}, A_{2, 1}, \dots, A_{2, d})^\top$, note that
\begin{equation*}
g^{-1}\left(\gamma_0 + \gamma_1\sum_{i = 1}^d \frac{A_{1, i}}{2^i} +\gamma_2\sum_{i = 1}^d \frac{A_{2, i}}{2^i}\right) = \bbeta^\top \bA^\otimes
\end{equation*}
for some BELIEF slopes $\bbeta$ by Theorem \ref{thm: 1-1GLM}. 

In general, $\bbeta$ has nonzero coefficients corresponding to interactions between bits of $U_1$ and $U_2$. For a concrete example, take $d = 1$, $\gamma_0 = 5$, and $\gamma_1 = \gamma_2 = 3$ with the logit link. In this case, we find the coefficient for the interaction between $A_{1,1}$ and $A_{2, 1}$ to be approximately $-0.18$. This value changes for different choices of $\gamma_0$, $\gamma_1$, and $\gamma_2$ and is zero when $\gamma_0 = 0, \gamma_1 = \gamma_2$. Moreover, due to Theorem \ref{thm: slopesAsDepthChanges}, BELIEF coefficients here are constant in $d$, meaning that these interactions do not vanish as the resolution grows. 

Whenever BELIEF interactions appear between bits of $U_1$ and $U_2$, the contribution of $U_2$ to $\Eb[B|U_1, U_2]$ must necessarily depend on the value of $U_1$. Intuitively, the existence of these interactions is apparent from the fact that $g$ has nonconstant derivative. More than simply affirming this fact, bitwise analysis produces rigorous statements on the magnitude of interaction at each resolution level. 

In summary, GLMs can impose interaction effects even when the structure purports to model main effects only, and vice versa. Moreover, as a consequence of Theorem \ref{thm: interactionNotEquiv}, logistic models regard interaction differently than, say, probit models. These properties hinder transparency in the sense of \cite{gelman2017beyond}, as the relationship between the choice of link and the interpretation of interaction effects is non-obvious. By considering bits up to a user-selected resolution, BELIEF may serve as a convenient linear benchmark that frames interaction directly in terms of the conditional mean.

\begin{remark}
On the subject of GLMs, the conditional independence properties of BELIEF relate to the notion of generating classes in log-linear models \citep{agresti2011categorical}. Generating classes indicate possible factorizations of the likelihood function for a contingency table and may be immediately translated into a graphical model. In contrast, BELIEF reflects a regression perspective, connecting conditional independence to coefficient values, meaning BELIEF provides direct information on the functional form of dependence. By way of the group operation introduced in Theorem \ref{thm: subsigma}, BELIEF accounts for nontrivial forms of conditional independence for multi-way interactions. Consequently, BELIEF provides alternative insights on conditional independence for binary interactions. We plan to investigate characterizations of conditional independence and graphical structure through BELIEF slopes, as well as their connections to log-linear models and other forms of GLMs.
\end{remark} 

\begin{remark}
In our exploration of binary dependence, we represent each ``binary'' bit $A$ as $\{-1, 1\}$-valued. For practical applications, binary variables are commonly encoded as $\tilde{A}\in \{0, 1\}$, where $\tilde{A} = (1+A)/2$. While this correspondence is trivial, the choice of representation affects the meaning of interaction terms. The product $A_1A_2$ equals one exactly when its factors are equal, while $\tilde{A}_1\tilde{A}_2$ equals one only if both its factors are one. (The distinction is analogous to the logic gates XNOR and AND.)

In a saturated scheme, both conventions yield equivalent probability models. For our theoretical purposes, the $\{-1, 1\}$ convention eliminates unnecessary asymmetry. As just one example, when using least squares estimation, all individual slopes necessarily have the same limiting variance (Theorem \ref{thm:asymnorm}), while $\{0, 1\}$ setup does not share this property.

As a final note, the ``hidden interaction'' principle of GLMs exists regardless of the binary encoding scheme. For instance, as we explain in this subsection, $g^{-1}(5 + (3/2)A_{1,1} + (3/2)A_{2,1}) = g^{-1}(2 + 3\tilde{A}_{1,1} + 3\tilde{A}_{2,1})$ has a nonzero BELIEF slope corresponding to the $A_{1,1}A_{2,1}$ interaction. 
\end{remark}

\subsection{The Case of Perfect Separation in Binary GLMs}\label{subsec: separation}

For a fixed binary sample $\{(\bA_i, B_i)\}_{i = 1}^n$, suppose that some function $f:\{-1, 1\}^p \rightarrow \{-1, 1\}$ satisfies $B_i = f(\bA_i)$ for all $1 \leq i \leq n$. In this case, the data are said to exhibit perfect separation, as the function $f$ explains all observed variability in the response. When this happens, GLMs for a binary response will fail in the sense that the MLE does not exist \citep{albert1984existence}, and with some non-informative prior, the posterior will not be proper either \citep{gelman1995bayesian,speckman2009existence}. This causes multiple problems, from conceptual to numerical and to practical \citep{rainey2016dealing}. It is counterintuitive that when there is a perfect predictor, which means that discriminability is best, maximum likelihood estimation fails to yield any slopes at all. Worse, when the data are close to perfect separation, the MLE is plagued by numerical instability. In such a situation, the GLM approach does not tell us which predictors yield perfect or near-perfect separation of the response. It is not apparent whether the failure is due to a lack of information or, rather, an excess via perfect separation.

While GLMs struggle with the above issues in perfect separation, BELIEF is able to avoid these problems and recover the relationship between $\bA_i$ and $B_i$ in the data. In Theorem~\ref{thm:sigma} and Theorem~\ref{thm:bounds}, we state how the least squares slope $\bbeta$ and its estimate $\hat{\bbeta}$ satisfy both $\norm{\Hb_{2^p}\bbeta}_\infty \leq 1$ and $\norm{\Hb_{2^p}\hat{\bbeta}}_\infty \leq 1$, which characterize a hypercube in $\mathbb{R}^{2^p}$ centered at $\zero$. As the next two results explain, the geometry of these constraints relate to the phenomenon of perfect separation in the data.

We first provide the result below on distributional perfect separation.

\begin{theorem}
\label{thm:hypercube_dist}
Let $\{(\bA, B)\}$ be binary variables where $\Eb[\bA^{\otimes}\bA^{\otimes^\top}]$ is positive definite. Then the following are equivalent:
\begin{enumerate}
    \item $\norm{\bbeta}_2 = 1.$
    \item $\bbeta$ is a vertex of the hypercube $\{\bx\in\mathbb{R}^{2^p}:\norm{\Hb_{2^p}\bx}_\infty \leq 1\}.$
    \item $B =\Eb[B|\bA]=\bbeta^\top\bA^{\otimes}$. In particular, $\Var(B)=\Var(\bbeta^\top\bA^{\otimes})=\Var(\Eb[B|\bA]).$ 
\end{enumerate}
\end{theorem}

Similarly, we have a version of the results for perfect separation in data.
\begin{theorem}
\label{thm:hypercube}
Let $\{(\bA_i, B_i)\}_{i = 1}^n$ be a binary sample where $\mathcal{A}^\top\mathcal{A}$ is positive definite. Then the following are equivalent:
\begin{enumerate}
    \item $\norm{\hat{\bbeta}}_2 = 1$
    \item $\hat{\bbeta}$ is a vertex of the hypercube $\{\bx\in\mathbb{R}^{2^p}:\norm{\Hb_{2^p}\bx}_\infty \leq 1\}$
    \item $B_i = \hat{\bbeta}^\top\bA_i^\otimes$ for all $1 \leq i \leq n$.
\end{enumerate}
\end{theorem}

Theorems \ref{thm:hypercube_dist} and \ref{thm:hypercube} show that in the situation of perfect separation, either theoretically or empirically, BELIEF slopes are well-defined and unique as parameters and as estimates. Moreover, note that $\Var(B)=\Var(\Eb[B|A])+\Eb[\Var(B|A)]$. Therefore, under perfect separation, BELIEF explains all the variance in the response, as a sensible model is expected to do in such a case. In this way, BELIEF overcomes the statistical and computational issues in binary GLM.

The geometric connections are also meaningful. There are $2^{2^p}$ vertices of the hypercube, and there are $2^{2^p}$ events in the $\sigma$-field generated by $\bA$ by Lemma~\ref{lem: pd}. Whenever $B = \Eb[B|\bA] = \bbeta^\top\bA^\otimes$, the variable $B$ can be thought of as a $-1/1$ indicator function of some event in $\sigma(\bA)$, where each distinct event corresponds to a different slope $\bbeta$. In light of Theorem \ref{thm:hypercube_dist}, these slopes are precisely the vertices of the hypercube $\{\bx\in\mathbb{R}^{2^p}:\norm{\Hb_{2^p}\bx}_\infty \leq 1\}.$ Therefore, BELIEF slope vectors provide clear description of this flavor of dependency.

\subsection{Approximating Arbitrary Dependency with BELIEF}
\label{subsec: doublelimit}
The BELIEF framework provides a promising strategy to approximate arbitrary dependency in the data. The main theorem of this section reassures us that decomposition into binary variables is a powerful tool for understanding dependence. Consider the simple case with $U$ and $V$ both supported within $[-1,1]$. With the binary expansion in Lemma~\ref{lem: binary:approx} up to depth $D$, we can decompose \textit{any} complex dependency between $U$ and $V$ into linear forms of dependency between each binary interaction in $\sigma$-fields generated from the binary expansions.

Specifically, let $U, V$ be random variables in $[-1, 1]$, not necessarily uniform. Let $\{A_d\}_{d = 1}^\infty$ and $\{B_d\}_{d = 1}^\infty$ be binary variables in $\{-1, 1\}$ such that
\begin{align*}
U = \sum_{d = 1}^\infty \frac{A_d}{2^d} \quad {\rm and}\quad
V = \sum_{d = 1}^\infty \frac{B_d}{2^d}.
\end{align*}
To ensure unique binary representation, we require that the binary expansion of rational numbers in $(-1, 1]$ end in repeating ones. For example, without this stipulation, the event $\{U = 0\}$ corresponds to both $\{A_1 = 1, A_d = -1 \text{ for } d > 1\}$ and $\{A_1 = -1, A_d = 1 \text{ for } d > 1\}$. Furthermore, let $\mathcal{F}_D = \sigma(\{A_d: 1 \leq d \leq D\})$ be the $\sigma$-field representing information in $U$ up to binary expansion depth $D$.  Hence $\{\mathcal{F}_j, j= 1, 2, \ldots \}$ forms an \textit{information filtration} for multi-resolution inference, where $j$ indexes the resolution level, as defined in the general framework by \cite{li2021multi}.  

\begin{theorem}
The following double limit holds for any $U$ and $V$ supported within $[-1,1]:$
\label{thm:doublelimit}
\begin{align*}
\Eb[V|U] \stackrel{a.s.}{=} \lim_{(D_1, D_2)\rightarrow\infty}\Eb\left[\sum_{d = 1}^{D_2} \frac{B_d}{2^d}\Bigg|\mathcal{F}_{D_1}\right] ,
\end{align*}
where $D_1$ denotes the binary expansion depth of $U$, and $D_2$ denotes the binary expansion depth of $V$.
\end{theorem}
\noindent
Theorem \ref{thm:doublelimit} assures us for predicting $V$ given $U$, using their binary expansions we can get arbitrarily close to using them directly. With the binary expansion in Lemma~\ref{lem: binary:approx} up to depth $D_1=D_2=D$, we can decompose \textit{any} complex dependency between $U$ and $V$ into linear forms of dependency in the BELIEF between each binary interaction in $\sigma(V_D)=\sigma(B_1,\ldots,B_D)$ on $\sigma(U_D)=\sigma(A_1,\ldots,A_D).$ In particular, conditional joint distributions of bits $B_1,\ldots,B_D$ are easily studied. This procedure can be represented graphically as in Figure~\ref{fig: ben} of Section~\ref{sec: discussion}.

So far, we have discussed the conditional expectation $\Eb[V|U]$ under the effect of binary expansion for both $V$ and $U$ simultaneously, showing uniform convergence of the approximation $\Eb\left[\sum_{d = 1}^{D_2} \frac{B_d}{2^d}\Big|\mathcal{F}_{D_1}\right]$ in $D_1$ and $D_2$. When applying BELIEF to probability models such as the GLM, we similarly replace continuous variables with their binary expansion in order to analyze bitwise slopes. Intuitively, the error of this approximation should tend to zero as the binary expansion depth increases.

For example, let $U_1, \dots, U_p \sim \text{Unif}(-1, 1)$, and let $A_{j,d}\in \{-1, 1\}$ denote the $d$th bit of $U_j$, where we treat the number of uniform predictors $p$ as a constant. Define the approximation $U_{j, D} = \sum_{d = 1}^D \frac{A_{j, d}}{2^d}$, where clearly $|U_j - U_{j, D}| = O_{\Pb}(2^{-D})$. In a typical application of BELIEF, we wish to approximate some function $g(U_1, \dots, U_p)$ via binary expansion of the arguments, such as in the modeling scheme $\Pb(B = 1| U_1, \dots, U_p) = g(U_1, \dots, U_p)$. If $g$ is a continuous function, then this approximation is correct in the sense that $g(U_{1, D}, \dots, U_{p, D}) \rightarrow g(U_1, \dots, U_p)$ almost surely as $D\rightarrow \infty$. Yet, given that the number of BELIEF slopes grows exponentially in $D$, we are also interested in the rate of this convergence, for fear that we pay more in parametrization than we receive in precision. It turns out that in many common applications such as GLMs, this rate of convergence is indeed exponential in $D$.

\begin{theorem}\label{thm: generalExponentialRate}
Let $U_1, \dots, U_p \sim \text{Unif}(-1, 1)$, let $A_{j,d}\in \{-1, 1\}$ denote the $d$th bit of $U_j$, and define $U_{j, D} = \sum_{d = 1}^D \frac{A_{j, d}}{2^d}$. Assume $g:(-1, 1)^p\rightarrow \mathbb{R}$ is a differentiable function. Then $|g(U_1, \dots, U_p) - g(U_{1, D}, \dots, U_{p, D})| = O_{\Pb}(2^{-D})$.
\end{theorem}

Importantly, the above result applies to functions of continuous, non-uniform random variables by way of the inverse CDF transformation. Additional assumptions on $g$ can yield stronger statements on the rate, such as a uniform rate in the case that $g$ is Lipschitz.

\section{Practical Considerations and Improvements}
\label{sec: improvementAndPractical}

\subsection{Dimensionality of BELIEF}\label{subsec: scalability}
One important practical concern is on the number of slopes in BELIEF. With $p$ covariates and each expanded to the $D$-th bit, we have $2^{pD}$ interactions in $\bA^\otimes.$ Thus, the number of BELIEF slopes grows exponentially in the number of covariates and bits, and the number of slopes may easily exceed the sample size. As with the conceptually similar $k$-nearest neighbor algorithm, a low sampling density in high dimensions leads to diminishing predictive performance \citep{htf2009}. A lack of structural assumptions is a liability in high dimensions, when the model asserts a level of granularity that the sample cannot hope to inform. Moreover, performing least squares with $2^{pD}$ predictors requires a computational complexity of $O(2^{3pD})$ in this current form. Another important issue is information loss from the binary expansion approximation to the joint distribution in Lemma~\ref{lem: binary:approx} as $p$ grows, whose $\ell^2$ norm is $O(\sqrt{p}2^{-D})$. These dimensionality problems are a faithful reflection of the hardness of the problem: to study arbitrary dependence, in growing dimensions, with no free lunch. These issues are also similar to the overparametrization problem in deep neural networks \citep{goodfellow2016deep}. However, strategies are needed to accommodate the dimensionality in BELIEF.

For one, the dimensionality issue can be combated through existing regularization methods, building on the contents of Section~\ref{subsec: regularization}. One might consider other dimension reduction methods such as screening \citep{fanlv2008} or resampling of randomly selected variables, as in the random forest algorithm \citep{htf2009}. A tempting third option involves fitting BELIEF slopes subject to structural constraints, before iteratively relaxing these constraints and comparing model selection criteria. Other approaches include borrowing theory and methods from wavelet analysis \citep{antoniadis2012wavelets} or tensor product models \citep{Wahba1994}. We will study dimension reduction in BELIEF carefully in a future paper, particularly structural assumptions such as sparsity \citep{sur2019modern}.

\subsection{Improvement with Regularization}\label{subsec: regularization}
One effective way to improve the prediction of regression models is to consider regularization or placing prior distributions over the slopes. In this section, we consider the use of shrinkage for estimating $\bbeta$. Recall that Theorem \ref{thm:condind} explains how zero-valued slopes correspond to conditional independence given subcollections of bits. Thus, shrinking some slopes towards zero can be thought of as moving closer to imposing various conditional independence assumptions between $\bA$ and $B$. 

Unlike shrinkage estimations with unbounded parameters, here we need to be mindful to ensure that we do not end up with illegitimate estimates of probabilities that go outside of the unit interval. For example, consider the ridge regression estimator
\begin{align}\label{eq:ridgr}
\hat{\bbeta}_{rid, \lambda} := \argmin_{\bbeta} \norm{\mathcal{B} - \mathcal{A}\bbeta}_2^2 + \lambda\norm{\bbeta}_2^2
\end{align}
for some $\lambda > 0$. The following result ensures that the minimization in $\eqref{eq:ridgr}$ always produces valid BELIEF slopes. 
\begin{theorem}
\label{thm:ridge}
The estimator $\hat{\bbeta}_{rid, \lambda}$ satisfies $\norm{\Hb_{2^p}\hat{\bbeta}_{rid, \lambda}}_\infty \leq 1,$
or equivalently, for any row $\bv$ of $\Hb_{2^p}$, 
\begin{align*}
\hat{\bP}(B = r|\bA^{\otimes} = \bv) & = \frac{1 + r\hat{\bbeta}_{rid}^\top\bv}{2} \in [0, 1], \quad {\rm for}\  r=\pm 1.
\end{align*}
Furthermore, whenever there is no observation $\bA^{\otimes}_i$ such that $\bA^{\otimes}_i = \bv$,  
$
\hat{\Pb}(B = r|\bA^{\otimes} = \bv) = 0.5 , \text{ for }  r=\pm 1.
$\end{theorem}
\noindent
Hence, the condition $\mathcal{A}^\top\mathcal{A}>0$ is unnecessary for the existence of $\hat{\bbeta}_{rid, \lambda}$, and the tuning parameter $\lambda$ reflects how much we wish to bias $\hat{\bbeta}_{rid, \lambda}$ towards zero. 

{One might also consider an analogous estimator $\hat{\bbeta}_{LASSO,\lambda}$ based on the LASSO penalty \citep{tibshirani96}. While we are unaware of any tractable general closed form representation of LASSO estimators, any function of binary inputs can be written as a multilinear function of those inputs, in the sense of Lemma~\ref{lem: linearity}. We emphasize a model-centric viewpoint of this fact, though another interesting consequence is that $\hat{\bbeta}_{MP}, \hat{\bbeta}_{rid, \lambda},$ $\hat{\bbeta}_{LASSO,\lambda}$, and any regularized estimator under BELIEF all admit some closed, multilinear representation in their binary data inputs. For the regularized estimators in particular, we wonder whether the multilinear representation could inspire an algorithm for computing them that, by taking advantage of bitwise operations, is faster than those commonly used in the case of generic data.}

\subsection{Other Practical Considerations}

One practical concern is the choice of the expansion depth $D$. In the simulations in \cite{zhang2021beauty}, it is found that with $D=3$, the test based on the binary expansion has a higher power than the linear model based tests for Gaussian data. Moreover, data studies in \cite{xiang2022pairwise} show that using $D=2$ can already provide many interesting findings. A general optimal choice of $D$ for some specific dependency should come from a trade-off between them and $n,p,$ and the signal strength. This would be an important problem for future studies.

Our perspective on continuous or other discrete response variables is similar to that of predictors: there is no conceptual barrier to fitting a separate BELIEF model for each bit in the binary expansion of the response. The information from these models can be aggregated to form a generic modeling strategy. The primary obstacle is instead computational, in the sense of the large number of parameters that would need to be estimated. This problem is worsened if the object of interest is the conditional distribution of the response, as opposed to just its conditional expectation, as an additional BELIEF model would be required for each interaction of bits from the response. While beyond our current scope, we believe a principled regularization of BELIEF slopes to be the solution to this type of issue.

We would also like to remark on the computational expense of fitting BELIEF. Because the bit structure is at the center of current computing systems, and because data are stored in computers as data bits, the steps of binary expansion and creation of $\bA_i^\otimes$ can be performed efficiently. As \cite{zhang2019bet} and \cite{zhao2021fast} show, algorithms based on binary expansion statistics can outperform existing nonparametric methods in speed by orders of magnitude.

\section{Empirical Investigations}\label{sec: empirical}

\subsection{Simulation Studies}
In this section, we study the performance of BELIEF in classification problems by comparing the receiver operating characteristic (ROC) curve with logistic regression, support vector machine (SVM) and random forest in three simulated examples with two continuous covariates $X_1$ and $X_2$ and a binary response $B$. For BELIEF, we consider the binary expansion with empirical CDF of marginal distributions up to depth 3 (64 slope parameters). For comparison, we consider a random forest with default settings and a logistic regression with the interaction between the two covariates (3 slope parameters).

The sample size of the simulations is set to be {12,288}, where {8,192} are randomly selected to form the training data, and the remaining {4,096} are used as test data to evaluate the classification accuracy. The data are simulated as i.i.d. observations from the following three scenarios. 
\begin{enumerate}
    \item Independent covariates and a linear predictor:\\ $X_1,X_2 \stackrel{i.i.d.}{\sim}\text{Unif}(-1,1)$ and $\logit(\Pb(B=1|X_1,X_2))=2X_1+X_2.$
    \item Independent covariates and a quadratic predictor:\\ $X_1,X_2 \stackrel{i.i.d.}{\sim}\cN(0,1)$ and $\logit(\Pb(B=1|X_1,X_2))=X_1^2+X_2^2.$
    \item Circularly dependent covariates and a cosine predictor:\\ $\theta\sim \text{Unif}[-\pi,\pi], \epsilon_1,\epsilon_2 \stackrel{i.i.d.}{\sim} \cN(0,1), X_1=\cos{\theta}+0.2\epsilon_1,X_2=\sin{\theta}+0.2\epsilon_2$ and $\logit(\Pb(B=1|X_1,X_2))=3\cos{(\pi(X_1+X_2))} .$
\end{enumerate}

\begin{figure}[ht]
\begin{center}
	\includegraphics[width=400pt]{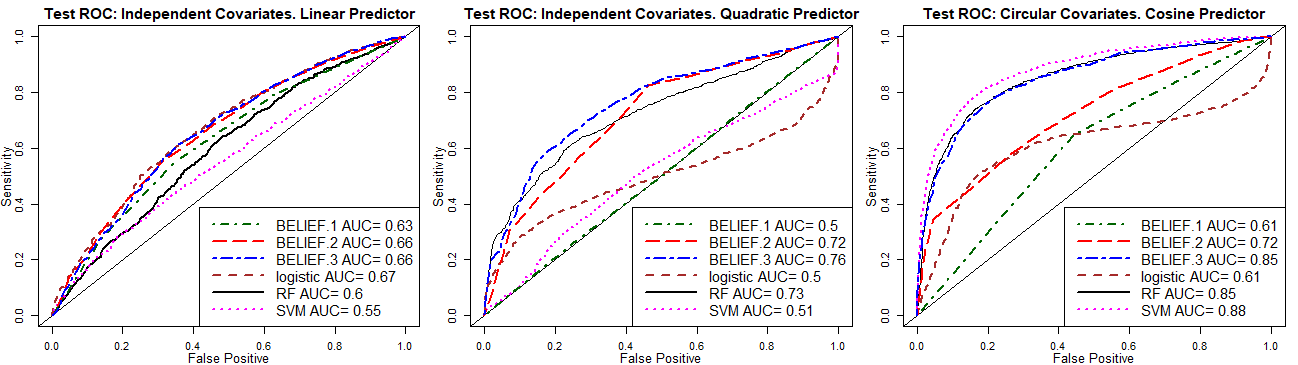}
\end{center}
{\footnotesize\caption{
		\footnotesize
		Comparison of ROC curves of different methods. BELIEF shows robust performance in all cases.}\label{fig: simu}
}
\end{figure}

The resulting numerical comparison of ROC curves is summarized in Figure~\ref{fig: simu}. In the left panel, the logistic regression shows the best ROC curve and has the highest area under curve (AUC). This is expected, since the link function in the logistic regression is well-specified and the predictor is indeed linear. The ROC curve of BELIEF at each depth shows a competitive performance in this case. This fact demonstrates the ability of BELIEF to approximate the optimal model.

\begin{table}
\caption{\footnotesize {Significant slopes fitted at depth $2$ in example 1 with $99\%$ confidence intervals with Bonferroni correction.}}
\label{tab: simu_e1d2}
\begin{center}
\begin{tabular}{r| c c }
Significant Binary Interaction & Slope & Confidence Interval \\ 
 \hline
 1 & 0.599 & (0.570, 0.628)\\
 $A_{2,2}$ & 0.031 & (0.002,  0.061)\\
 $A_{2,1}$ & 0.067 & (0.038, 0.096) \\
 $A_{1,2}$ & 0.069 & (0.039, 0.098)\\
 $A_{1,1}$ & $0.160$ & $(0.130, 0.189)$
\end{tabular}
\end{center}
\end{table}

Compared to other nonparametric methods, the inference of slopes in BELIEF provides clear interpretation on the form of the dependency. {Using a cutoff of $\alpha = 0.01$ and Bonferroni correction,} Table~\ref{tab: simu_e1d2} shows that, {out of 16 total interactions}, the first two bits of {both $X_1$ and $X_2$ have significantly positive main effects on the response. This corresponds to the monotone relationship between the predictors $X_1$, $X_2$ and $\Pb(B=1|X_1,X_2)$.}

In the middle panel of Figure \ref{fig: simu}, since the predictor is nonlinear, the logistic regression has a much worse performance despite the correct link function. Nonparametric methods, such as random forest, SVM and BELIEF, show a much more robust performance instead. Although BELIEF at depth 1 is not capable of capturing nonlinear dependency, BELIEF at depth 2 and 3 both provide a good ROC curve with high AUC.  

The inference of the slopes provides clear interpretation of the nonlinear relationship. At depth $2$, three slopes {out of 16} in BELIEF are significant {at the $\alpha = 0.01$ level after Bonferroni correction}, as shown in Table~\ref{tab: simu_e2d2}. Note the interaction $A_{1,1}A_{1,2}$ corresponds to contrasting the observations in the inner half of the distribution of $X_1$ versus the outer half. This interaction thus indicates a quadratic functional relationship between $X_1$ and the response $B$. The significance of the slope of $A_{2,1}A_{2,2}$ hints the same story. In addition, the interaction  $A_{1,1}A_{1,2}A_{2,1}A_{2,2}$ suggests a circular contour for the conditional probability. These slopes thus clearly describe the dependency between $(X_1,X_2)$ and $B.$

The significant BELIEF slopes also show the possibility of conditional independence. From results in Section~\ref{subsec: condind}, we see that
$1, A_{1,1}A_{1,2}, A_{2,1}A_{2,2}, A_{1,1}A_{1,2}A_{2,1}A_{2,2}$ form a proper subgroup of the multiplicative group of all binary interactions. Therefore, $B$ is conditionally independent of any binary interaction given the interactions of the first two bits from each predictor, i.e., $B \indep \bA^{\otimes}|(A_{1,1}A_{1,2},A_{2,1}A_{2,2}).$ In practice, this means all dependency between the response $B$ and the predictors are within $(A_{1,1}A_{1,2},A_{2,1}A_{2,2})$ which approximates the circular form of dependency.

\begin{table}
\caption{\footnotesize {Significant slopes fitted at depth $2$ in example 2 with $99\%$ confidence intervals with Bonferroni correction.}}
\label{tab: simu_e2d2}
\begin{center}
\begin{tabular}{r| c c }
Significant Binary Interaction & Slope & Confidence Interval \\ 
 \hline
 1 & 0.568 & (0.539, 0.598)\\
 $A_{1,1}A_{1,2}$ & 0.187 & (0.157,0.216)\\
 $A_{2,1}A_{2,2}$ & 0.186 & (0.156,0.215)\\
 $A_{1,1}A_{1,2}A_{2,1}A_{2,2}$ & -0.073 & (-0.102,-0.043)
\end{tabular}
\end{center}
\end{table}

In example 3, the logistic regression is again incapable of modeling the nonlinear relationship despite the correct link. For BELIEF, the complex dependency in the data requires a higher depth in the binary expansion to model the dependency. BELIEF with depth 3 shows a higher ROC curve than other methods, as shown in the right panel of Figure \ref{fig: simu}.

In all of the three examples above, BELIEF shows a competitive performance among existing classification methods. Compared to parametric methods, BELIEF is more robust against nonlinear forms of dependency through a distribution-free approach. Compared to nonparametric methods, BELIEF provides meaningful interpretations and useful inference of the forms of dependency through slopes and subgroups of binary interactions. These features shows the potential of BELIEF as a useful alternative method in practice.

{

\subsection{Data Example}

In Section \ref{subsec: illustrative}, we allude to a longstanding discussion in the political science literature regarding the necessity and use of interaction terms in logistic regression \citep{berry2010testing, rainey2016compression}. While we give theoretically-oriented comments on the matter in Section \ref{subsec: generalresultsforglms}, we take this opportunity to speak practically, borrowing the logistic example from \cite{rainey2016compression} and comparing with BELIEF.

The data originally come from \cite{russett2001triangulating}, who used a logistic model to predict conflict between countries. Conflict between two countries in a given year constitutes the binary response. Binary predictors include whether the countries were previously allies and whether the countries are contiguous. Other predictors include the log ratio of capabilities between the two countries, the log distance between capitals or major ports, and the smaller of the two countries' polity scores, which measures the balance of autocracy versus democracy. Compared to the original analysis of this data, we omit the binary control describing whether both countries are minor powers due to redundancy: by a separate logistic regression, we find that this variable is very well-predicted by the other variables already included.

We re-analyze the data with BELIEF with depth $d=1,$ i.e., we consider the two binary predictors described above and the first bit of the three continuous variables. The inference of BELIEF slopes is based on the Bonferroni correction at $0.01$ over $2^5=32$ slopes. According to this adjustment, all five main effects are significant, and six two-way binary interactions are significant. All other slopes are not Bonferroni significant. We list below the six two-way interactions that are significant:
		\begin{itemize}
			\item whether the countries were previously allies and whether the countries are contiguous;
			\item whether the countries were previously allies and the first bit of capability ratio;
			\item whether the countries are contiguous and the first bit of capability ratio;
			\item whether the countries were previously allies and the first bit of  the smaller of the two countries' polity scores;
			\item whether the countries are contiguous and the first bit of the distance between capitals;
			\item the first bit of capability ratio and the first bit of the smaller of the two countries' polity scores.
		\end{itemize}
Although our analysis is only based on $d=1$, the inference above has already provided interesting and meaningful interpretations. All significant two-way interactions seem to make political scientific sense. Moreover, the interaction between the first bits of two continuous predictors indicate a nonlinear form of dependency. We would like to emphasize that the above meaningful inference is achieved without making any assumptions about the relationship or the choice of a link function. Therefore, inference with BELIEF could serve as a useful alternative to GLM, offering new insights.

}

\section{Further Explorations}\label{sec: discussion}

\begin{figure}[ht]
	\begin{center}
		\includegraphics[width=220pt]{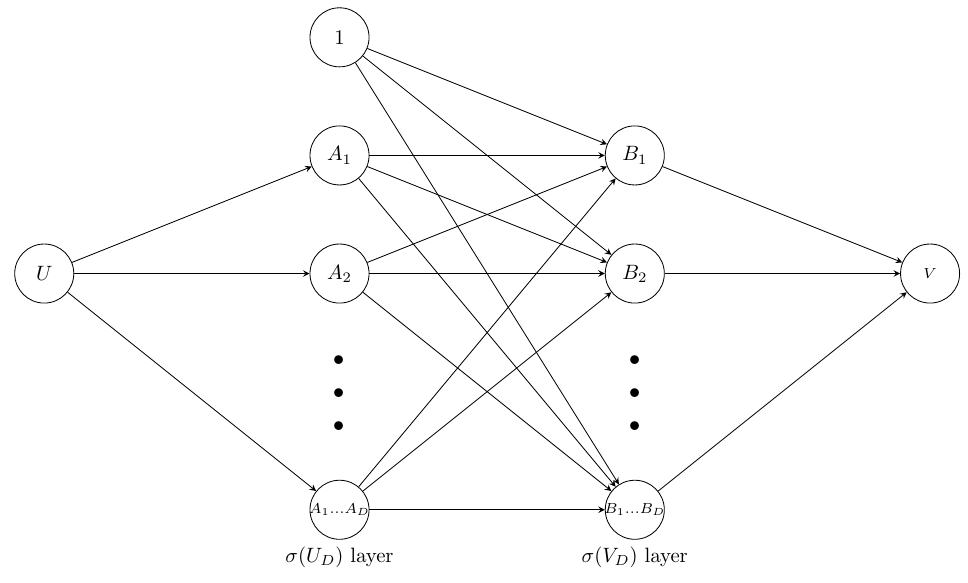}
	\end{center}
	{\footnotesize\caption{
			\footnotesize
			The network representation of the BELIEF modeling of dependency between $U$ and $V$.}\label{fig: ben}
	}
\end{figure}

BELIEF restates the dependence of binary variables in the familiar language of linear models, providing an interpretable parameterization of dependence up to a predetermined resolution. While BELIEF is formulated using the structure of binary variables, one can apply BELIEF to continuous variables even when the marginal distribution is unknown. We furnish a host of theoretical results showing the convenience of binary representation: the unintuitive properties of GLMs are of particular interest. Further, we prototype estimation with BELIEF using real-world data, and we articulate interesting properties of binary variables, including comparisons with the Gaussian family.

As a promising future direction, we wish to remark on the connections between BELIEF and deep learning. Figure~\ref{fig: ben} illustrates the interesting structural similarity between dependence with BELIEF and artificial neural networks \citep{goodfellow2016deep}. The input layer of this network is formed by observations in $U$. The nodes in the $\sigma(U_D)$ layer decomposes the low-resolution information in the distribution of $U$ into those of the binary interactions of data bits $A_1,\ldots,A_1A_2,\ldots,A_1\cdots A_D.$ The $\sigma(V_D)$ layer fits the BELIEF of each of $B_1,\ldots,B_1B_2,\ldots,B_1\cdots B_D$ on $\sigma(U_D)$ to obtain their conditional expectations, where weights of nodes are the slopes. These conditional expectations can then be aggregated to provide distributional statements of the output layer of $V$. These layers constitute an object one might call a binary expansion network (BEN). In addition, while being mindful that association does not imply causation, Theorem \ref{thm:condind} indicates the close relationship between BELIEF and Gaussian properties as applied to causal inference, where conditional dependence manifests conspicuously in BELIEF slopes. An increased emphasis in recent research of neural networks is about their interpretability \citep{murdoch2019definitions,rudin2019stop}. Motivated by the simple and easily interpretable forms of the BELIEF, we propose exploring the connections to general neural networks in the future to improve the interpretability of deep learning, as well as other possible improvements. 

As in deep neural networks, one major challenge in the BELIEF framework is the reduction of dimensions, as discussed in Section~\ref{subsec: scalability}. Existing and novel theory and methods on dimension reduction could play essential roles in further studies under this framework. We welcome discussions and suggestions on this topic for deeper understanding of dependence and useful procedures in practice.

\begin{acks}[Acknowledgments]
The authors thank Alan Agresti, Tony Cai, Edgar Dobriban, Chao Gao, Fang Han, Lihua Lei, Tengyuan Liang, Peter McCullagh, and Xianyang Zhang for their invaluable comments and suggestions. The proof of Theorem~3.8 is based on techniques developed in personal communications with Xianyang Zhang. In addition, the authors thank Wan Zhang for her helpful mathematical and computational insights. The authors further thank the AOS reviewers for their constructive comments.
\end{acks}

\begin{funding}
Brown's and Zhang’s research is partially supported by NSF grants DMS-1916237 and DMS-2152289. Meng thanks NSF for partial financial support.
\end{funding}

\begin{supplement}
\stitle{Additional Theory and Proofs}
\sdescription{Further results on binary variables, as well as proofs of theorems, lemmas, and corollaries.}
\end{supplement}


\bibliographystyle{imsart-number} 
\bibliography{bibliography}       

\begin{thebibliography}{54}

\bibitem{agresti2011categorical}
\begin{bbook}[author]
\bauthor{\bsnm{Agresti},~\bfnm{Alan}\binits{A.}} \AND \bauthor{\bsnm{Kateri},~\bfnm{Maria}\binits{M.}}
(\byear{2011}).
\btitle{Categorical data analysis}.
\bpublisher{Springer}.
\end{bbook}
\endbibitem

\bibitem{albert1984existence}
\begin{barticle}[author]
\bauthor{\bsnm{Albert},~\bfnm{Adelin}\binits{A.}} \AND \bauthor{\bsnm{Anderson},~\bfnm{John~A}\binits{J.~A.}}
(\byear{1984}).
\btitle{On the existence of maximum likelihood estimates in logistic regression models}.
\bjournal{Biometrika}
\bvolume{71}
\bpages{1--10}.
\end{barticle}
\endbibitem

\bibitem{angrist2009mostly}
\begin{bbook}[author]
\bauthor{\bsnm{Angrist},~\bfnm{Joshua~D}\binits{J.~D.}} \AND \bauthor{\bsnm{Pischke},~\bfnm{J{\"o}rn-Steffen}\binits{J.-S.}}
(\byear{2009}).
\btitle{Mostly harmless econometrics: An empiricist's companion}.
\bpublisher{Princeton university press}.
\end{bbook}
\endbibitem

\bibitem{antoniadis2012wavelets}
\begin{bbook}[author]
\bauthor{\bsnm{Antoniadis},~\bfnm{Anestis}\binits{A.}} \AND \bauthor{\bsnm{Oppenheim},~\bfnm{Georges}\binits{G.}}
(\byear{2012}).
\btitle{Wavelets and statistics}.
\bpublisher{Springer Science \& Business Media}.
\end{bbook}
\endbibitem

\bibitem{bahadur1961representation}
\begin{barticle}[author]
\bauthor{\bsnm{Bahadur},~\bfnm{Raghu~Raj}\binits{R.~R.}}
(\byear{1961}).
\btitle{A representation of the joint distribution of responses to n dichotomous items}.
\bjournal{Studies in item analysis and prediction}
\bpages{158--168}.
\end{barticle}
\endbibitem

\bibitem{barber2020distribution}
\begin{barticle}[author]
\bauthor{\bsnm{Barber},~\bfnm{Rina~Foygel}\binits{R.~F.}}
(\byear{2020}).
\btitle{Is distribution-free inference possible for binary regression?}
\bjournal{Electronic Journal of Statistics}
\bvolume{14}
\bpages{3487--3524}.
\end{barticle}
\endbibitem

\bibitem{foygel2021limits}
\begin{barticle}[author]
\bauthor{\bsnm{Barber},~\bfnm{R~F.}\binits{R.~F.}}, \bauthor{\bsnm{Candes},~\bfnm{Emmanuel~J}\binits{E.~J.}}, \bauthor{\bsnm{Ramdas},~\bfnm{Aaditya}\binits{A.}} \AND \bauthor{\bsnm{Tibshirani},~\bfnm{Ryan~J}\binits{R.~J.}}
(\byear{2021}).
\btitle{The limits of distribution-free conditional predictive inference}.
\bjournal{Information and Inference: A Journal of the IMA}
\bvolume{10}
\bpages{455--482}.
\end{barticle}
\endbibitem

\bibitem{berry2010testing}
\begin{barticle}[author]
\bauthor{\bsnm{Berry},~\bfnm{William~D}\binits{W.~D.}}, \bauthor{\bsnm{DeMeritt},~\bfnm{Jacqueline~HR}\binits{J.~H.}} \AND \bauthor{\bsnm{Esarey},~\bfnm{Justin}\binits{J.}}
(\byear{2010}).
\btitle{Testing for interaction in binary logit and probit models: is a product term essential?}
\bjournal{American Journal of Political Science}
\bvolume{54}
\bpages{248--266}.
\end{barticle}
\endbibitem

\bibitem{beliefsupp}
\begin{barticle}[author]
\bauthor{\bsnm{Brown},~\bfnm{Benjamin}\binits{B.}}, \bauthor{\bsnm{Zhang},~\bfnm{Kai}\binits{K.}} \AND \bauthor{\bsnm{Meng},~\bfnm{Xiao-Li}\binits{X.-L.}}
\btitle{Supplement to ``BELIEF in Dependence: Leveraging Atomic Linearity in Data Bits for Rethinking Generalized Linear Models''}.
\end{barticle}
\endbibitem

\bibitem{buhlmann2011}
\begin{bbook}[author]
\bauthor{\bsnm{B\"{u}hlmann},~\bfnm{Peter}\binits{P.}} \AND \bauthor{\bparticle{van~de} \bsnm{Geer},~\bfnm{Sara}\binits{S.}}
(\byear{2011}).
\btitle{Statistics for high-dimensional data}.
\bpublisher{Springer}.
\end{bbook}
\endbibitem

\bibitem{buja2019}
\begin{barticle}[author]
\bauthor{\bsnm{Buja},~\bfnm{Andreas}\binits{A.}}, \bauthor{\bsnm{Brown},~\bfnm{Lawrence}\binits{L.}}, \bauthor{\bsnm{Berk},~\bfnm{Richard}\binits{R.}}, \bauthor{\bsnm{George},~\bfnm{Edward}\binits{E.}}, \bauthor{\bsnm{Pitkin},~\bfnm{Emil}\binits{E.}}, \bauthor{\bsnm{Traskin},~\bfnm{Mikhail}\binits{M.}}, \bauthor{\bsnm{Zhang},~\bfnm{Kai}\binits{K.}} \AND \bauthor{\bsnm{Zhao},~\bfnm{Linda}\binits{L.}}
(\byear{2019}).
\btitle{Models as Approximations {I}: Consequences Illustrated with Linear Regression}.
\bjournal{Statistical Science}
\bvolume{34}
\bpages{523--544}.
\bdoi{10.1214/18-STS693}
\end{barticle}
\endbibitem

\bibitem{chen2016fused}
\begin{barticle}[author]
\bauthor{\bsnm{Chen},~\bfnm{Yunxiao}\binits{Y.}}, \bauthor{\bsnm{Li},~\bfnm{Xiaoou}\binits{X.}}, \bauthor{\bsnm{Liu},~\bfnm{Jingchen}\binits{J.}} \AND \bauthor{\bsnm{Ying},~\bfnm{Zhiliang}\binits{Z.}}
(\byear{2016}).
\btitle{A fused latent and graphical model for multivariate binary data}.
\bjournal{arXiv preprint arXiv:1606.08925}.
\end{barticle}
\endbibitem

\bibitem{chernoff2004information}
\begin{barticle}[author]
\bauthor{\bsnm{Chernoff},~\bfnm{Herman}\binits{H.}}
(\byear{2004}).
\btitle{Information, for testing the equality of two probabilities, from the margins of the 2$\times$ 2 table}.
\bjournal{Journal of statistical planning and inference}
\bvolume{121}
\bpages{209--214}.
\end{barticle}
\endbibitem

\bibitem{chernozhukov2022high}
\begin{barticle}[author]
\bauthor{\bsnm{Chernozhukov},~\bfnm{Victor}\binits{V.}}, \bauthor{\bsnm{Chetverikov},~\bfnm{Denis}\binits{D.}}, \bauthor{\bsnm{Kato},~\bfnm{Kengo}\binits{K.}} \AND \bauthor{\bsnm{Koike},~\bfnm{Yuta}\binits{Y.}}
(\byear{2022}).
\btitle{High-dimensional Data Bootstrap}.
\bjournal{arXiv preprint arXiv:2205.09691}.
\end{barticle}
\endbibitem

\bibitem{cox1972analysis}
\begin{barticle}[author]
\bauthor{\bsnm{Cox},~\bfnm{David~R}\binits{D.~R.}}
(\byear{1972}).
\btitle{The analysis of multivariate binary data}.
\bjournal{Applied statistics}
\bpages{113--120}.
\end{barticle}
\endbibitem

\bibitem{dai2013multivariate}
\begin{barticle}[author]
\bauthor{\bsnm{Dai},~\bfnm{Bin}\binits{B.}}, \bauthor{\bsnm{Ding},~\bfnm{Shilin}\binits{S.}} \AND \bauthor{\bsnm{Wahba},~\bfnm{Grace}\binits{G.}}
(\byear{2013}).
\btitle{{Multivariate Bernoulli distribution}}.
\bjournal{Bernoulli}
\bvolume{19}
\bpages{1465 -- 1483}.
\bdoi{10.3150/12-BEJSP10}
\end{barticle}
\endbibitem

\bibitem{esary1967association}
\begin{barticle}[author]
\bauthor{\bsnm{Esary},~\bfnm{James~D}\binits{J.~D.}}, \bauthor{\bsnm{Proschan},~\bfnm{Frank}\binits{F.}} \AND \bauthor{\bsnm{Walkup},~\bfnm{David~W}\binits{D.~W.}}
(\byear{1967}).
\btitle{Association of random variables, with applications}.
\bjournal{The Annals of Mathematical Statistics}
\bvolume{38}
\bpages{1466--1474}.
\end{barticle}
\endbibitem

\bibitem{fanlv2008}
\begin{barticle}[author]
\bauthor{\bsnm{Fan},~\bfnm{Jianqing}\binits{J.}} \AND \bauthor{\bsnm{Lv},~\bfnm{Jinchi}\binits{J.}}
(\byear{2008}).
\btitle{Sure independence screening for ultrahigh dimensional feature space}.
\bjournal{Journal of the Royal Statistical Society: Series B}
\bvolume{70}
\bpages{849-911}.
\bdoi{10.1111/j.1467-9868.2008.00674.x}
\end{barticle}
\endbibitem

\bibitem{fan2015innovated}
\begin{barticle}[author]
\bauthor{\bsnm{Fan},~\bfnm{Yingying}\binits{Y.}}, \bauthor{\bsnm{Kong},~\bfnm{Yinfei}\binits{Y.}}, \bauthor{\bsnm{Li},~\bfnm{Daoji}\binits{D.}} \AND \bauthor{\bsnm{Zheng},~\bfnm{Zemin}\binits{Z.}}
(\byear{2015}).
\btitle{Innovated interaction screening for high-dimensional nonlinear classification}.
\bjournal{The Annals of Statistics}
\bvolume{43}
\bpages{1243--1272}.
\end{barticle}
\endbibitem

\bibitem{gelman1995bayesian}
\begin{bbook}[author]
\bauthor{\bsnm{Gelman},~\bfnm{Andrew}\binits{A.}}, \bauthor{\bsnm{Carlin},~\bfnm{John~B}\binits{J.~B.}}, \bauthor{\bsnm{Stern},~\bfnm{Hal~S}\binits{H.~S.}} \AND \bauthor{\bsnm{Rubin},~\bfnm{Donald~B}\binits{D.~B.}}
(\byear{1995}).
\btitle{Bayesian data analysis}.
\bpublisher{Chapman and Hall/CRC}.
\end{bbook}
\endbibitem

\bibitem{gelman2017beyond}
\begin{barticle}[author]
\bauthor{\bsnm{Gelman},~\bfnm{Andrew}\binits{A.}} \AND \bauthor{\bsnm{Hennig},~\bfnm{Christian}\binits{C.}}
(\byear{2017}).
\btitle{Beyond subjective and objective in statistics}.
\bjournal{Journal of the Royal Statistical Society: Series A (Statistics in Society)}
\bvolume{180}
\bpages{967--1033}.
\end{barticle}
\endbibitem

\bibitem{goodfellow2016deep}
\begin{bbook}[author]
\bauthor{\bsnm{Goodfellow},~\bfnm{Ian}\binits{I.}}, \bauthor{\bsnm{Bengio},~\bfnm{Yoshua}\binits{Y.}} \AND \bauthor{\bsnm{Courville},~\bfnm{Aaron}\binits{A.}}
(\byear{2016}).
\btitle{Deep learning}.
\bpublisher{MIT Press}.
\end{bbook}
\endbibitem

\bibitem{gupta2020distribution}
\begin{barticle}[author]
\bauthor{\bsnm{Gupta},~\bfnm{Chirag}\binits{C.}}, \bauthor{\bsnm{Podkopaev},~\bfnm{Aleksandr}\binits{A.}} \AND \bauthor{\bsnm{Ramdas},~\bfnm{Aaditya}\binits{A.}}
(\byear{2020}).
\btitle{Distribution-free binary classification: prediction sets, confidence intervals and calibration}.
\bjournal{Advances in Neural Information Processing Systems}
\bvolume{33}
\bpages{3711--3723}.
\end{barticle}
\endbibitem

\bibitem{htf2009}
\begin{bbook}[author]
\bauthor{\bsnm{Hastie},~\bfnm{Trevor}\binits{T.}}, \bauthor{\bsnm{Tibshirani},~\bfnm{Robert}\binits{R.}} \AND \bauthor{\bsnm{Friedman},~\bfnm{Jerome}\binits{J.}}
(\byear{2009}).
\btitle{The Elements of Statistical Learning: Prediction, Inference and Data Mining},
\bedition{2nd} ed.
\bpublisher{Springer Verlag.}
\end{bbook}
\endbibitem

\bibitem{he2013quantile}
\begin{barticle}[author]
\bauthor{\bsnm{He},~\bfnm{Xuming}\binits{X.}}, \bauthor{\bsnm{Wang},~\bfnm{Lan}\binits{L.}} \AND \bauthor{\bsnm{Hong},~\bfnm{Hyokyoung~Grace}\binits{H.~G.}}
(\byear{2013}).
\btitle{Quantile-adaptive model-free variable screening for high-dimensional heterogeneous data}.
\bjournal{The Annals of Statistics}
\bvolume{41}
\bpages{342--369}.
\end{barticle}
\endbibitem

\bibitem{lauritzen1996graphical}
\begin{bbook}[author]
\bauthor{\bsnm{Lauritzen},~\bfnm{Steffen~L}\binits{S.~L.}}
(\byear{1996}).
\btitle{Graphical models}.
\bpublisher{Clarendon Press}.
\end{bbook}
\endbibitem

\bibitem{lei2018distribution}
\begin{barticle}[author]
\bauthor{\bsnm{Lei},~\bfnm{Jing}\binits{J.}}, \bauthor{\bsnm{G’Sell},~\bfnm{Max}\binits{M.}}, \bauthor{\bsnm{Rinaldo},~\bfnm{Alessandro}\binits{A.}}, \bauthor{\bsnm{Tibshirani},~\bfnm{Ryan~J}\binits{R.~J.}} \AND \bauthor{\bsnm{Wasserman},~\bfnm{Larry}\binits{L.}}
(\byear{2018}).
\btitle{Distribution-free predictive inference for regression}.
\bjournal{Journal of the American Statistical Association}
\bvolume{113}
\bpages{1094--1111}.
\end{barticle}
\endbibitem

\bibitem{li2020nonparametric}
\begin{barticle}[author]
\bauthor{\bsnm{Li},~\bfnm{Jialu}\binits{J.}}, \bauthor{\bsnm{Zhang},~\bfnm{Wan}\binits{W.}}, \bauthor{\bsnm{Wang},~\bfnm{Peiyao}\binits{P.}}, \bauthor{\bsnm{Li},~\bfnm{Qizhai}\binits{Q.}}, \bauthor{\bsnm{Zhang},~\bfnm{Kai}\binits{K.}} \AND \bauthor{\bsnm{Liu},~\bfnm{Yufeng}\binits{Y.}}
(\byear{2022}).
\btitle{Nonparametric prediction distribution from resolution-wise regression with heterogeneous data}.
\bjournal{Journal of Business and Economics Statistics}
\bvolume{to appear}.
\end{barticle}
\endbibitem

\bibitem{li2021multi}
\begin{barticle}[author]
\bauthor{\bsnm{Li},~\bfnm{Xinran}\binits{X.}} \AND \bauthor{\bsnm{Meng},~\bfnm{Xiao-Li}\binits{X.-L.}}
(\byear{2021}).
\btitle{A multi-resolution theory for approximating infinite-$p$-zero-$n$: Transitional inference, individualized predictions, and a world without bias-variance tradeoff}.
\bjournal{Journal of the American Statistical Association}
\bvolume{116}
\bpages{353--367}.
\end{barticle}
\endbibitem

\bibitem{lin2000}
\begin{barticle}[author]
\bauthor{\bsnm{Lin},~\bfnm{Yi}\binits{Y.}}
(\byear{2000}).
\btitle{{Tensor product space ANOVA models}}.
\bjournal{The Annals of Statistics}
\bvolume{28}
\bpages{734 -- 755}.
\bdoi{10.1214/aos/1015951996}
\end{barticle}
\endbibitem

\bibitem{mccullagh2000invariance}
\begin{barticle}[author]
\bauthor{\bsnm{McCullagh},~\bfnm{Peter}\binits{P.}}
(\byear{2000}).
\btitle{Invariance and factorial models}.
\bjournal{Journal of the Royal Statistical Society: Series B (Statistical Methodology)}
\bvolume{62}
\bpages{209--256}.
\end{barticle}
\endbibitem

\bibitem{mccullagh2019generalized}
\begin{bbook}[author]
\bauthor{\bsnm{McCullagh},~\bfnm{Peter}\binits{P.}} \AND \bauthor{\bsnm{Nelder},~\bfnm{John~A}\binits{J.~A.}}
(\byear{2019}).
\btitle{Generalized linear models}.
\bpublisher{Routledge}.
\end{bbook}
\endbibitem

\bibitem{murdoch2019definitions}
\begin{barticle}[author]
\bauthor{\bsnm{Murdoch},~\bfnm{W~James}\binits{W.~J.}}, \bauthor{\bsnm{Singh},~\bfnm{Chandan}\binits{C.}}, \bauthor{\bsnm{Kumbier},~\bfnm{Karl}\binits{K.}}, \bauthor{\bsnm{Abbasi-Asl},~\bfnm{Reza}\binits{R.}} \AND \bauthor{\bsnm{Yu},~\bfnm{Bin}\binits{B.}}
(\byear{2019}).
\btitle{Definitions, methods, and applications in interpretable machine learning}.
\bjournal{Proceedings of the National Academy of Sciences}
\bvolume{116}
\bpages{22071--22080}.
\end{barticle}
\endbibitem

\bibitem{o2014analysis}
\begin{bbook}[author]
\bauthor{\bsnm{O'Donnell},~\bfnm{Ryan}\binits{R.}}
(\byear{2014}).
\btitle{Analysis of boolean functions}.
\bpublisher{Cambridge University Press}.
\end{bbook}
\endbibitem

\bibitem{rainey2016compression}
\begin{barticle}[author]
\bauthor{\bsnm{Rainey},~\bfnm{Carlisle}\binits{C.}}
(\byear{2016}).
\btitle{Compression and conditional effects: A product term is essential when using logistic regression to test for interaction}.
\bjournal{Political Science Research and Methods}
\bvolume{4}
\bpages{621--639}.
\end{barticle}
\endbibitem

\bibitem{rainey2016dealing}
\begin{barticle}[author]
\bauthor{\bsnm{Rainey},~\bfnm{Carlisle}\binits{C.}}
(\byear{2016}).
\btitle{Dealing with separation in logistic regression models}.
\bjournal{Political Analysis}
\bvolume{24}
\bpages{339--355}.
\end{barticle}
\endbibitem

\bibitem{ravikumar2010high}
\begin{barticle}[author]
\bauthor{\bsnm{Ravikumar},~\bfnm{Pradeep}\binits{P.}}, \bauthor{\bsnm{Wainwright},~\bfnm{Martin~J}\binits{M.~J.}} \AND \bauthor{\bsnm{Lafferty},~\bfnm{John~D}\binits{J.~D.}}
(\byear{2010}).
\btitle{High-dimensional {I}sing model selection using $\ell_1$-regularized logistic regression}.
\bjournal{The Annals of Statistics}
\bvolume{38}
\bpages{1287--1319}.
\end{barticle}
\endbibitem

\bibitem{rosenbaum1983central}
\begin{barticle}[author]
\bauthor{\bsnm{Rosenbaum},~\bfnm{Paul~R}\binits{P.~R.}} \AND \bauthor{\bsnm{Rubin},~\bfnm{Donald~B}\binits{D.~B.}}
(\byear{1983}).
\btitle{The central role of the propensity score in observational studies for causal effects}.
\bjournal{Biometrika}
\bvolume{70}
\bpages{41--55}.
\end{barticle}
\endbibitem

\bibitem{rudin2019stop}
\begin{barticle}[author]
\bauthor{\bsnm{Rudin},~\bfnm{Cynthia}\binits{C.}}
(\byear{2019}).
\btitle{Stop explaining black box machine learning models for high stakes decisions and use interpretable models instead}.
\bjournal{Nature Machine Intelligence}
\bvolume{1}
\bpages{206--215}.
\end{barticle}
\endbibitem

\bibitem{russett2001triangulating}
\begin{bbook}[author]
\bauthor{\bsnm{Russett},~\bfnm{Bruce~M}\binits{B.~M.}} \AND \bauthor{\bsnm{Oneal},~\bfnm{John~R}\binits{J.~R.}}
(\byear{2001}).
\btitle{Triangulating peace: Democracy, interdependence, and international organizations}.
\bpublisher{W. W. Norton}, \baddress{New York}.
\end{bbook}
\endbibitem

\bibitem{speckman2009existence}
\begin{barticle}[author]
\bauthor{\bsnm{Speckman},~\bfnm{Paul~L}\binits{P.~L.}}, \bauthor{\bsnm{Lee},~\bfnm{Jaeyong}\binits{J.}} \AND \bauthor{\bsnm{Sun},~\bfnm{Dongchu}\binits{D.}}
(\byear{2009}).
\btitle{Existence of the {MLE} and propriety of posteriors for a general multinomial choice model}.
\bjournal{Statistica Sinica}
\bvolume{19}
\bpages{731-748}.
\end{barticle}
\endbibitem

\bibitem{sur2019modern}
\begin{barticle}[author]
\bauthor{\bsnm{Sur},~\bfnm{Pragya}\binits{P.}} \AND \bauthor{\bsnm{Cand{\`e}s},~\bfnm{Emmanuel~J}\binits{E.~J.}}
(\byear{2019}).
\btitle{A modern maximum-likelihood theory for high-dimensional logistic regression}.
\bjournal{Proceedings of the National Academy of Sciences}
\bvolume{116}
\bpages{14516--14525}.
\end{barticle}
\endbibitem

\bibitem{thanei2018xyz}
\begin{barticle}[author]
\bauthor{\bsnm{Thanei},~\bfnm{Gian-Andrea}\binits{G.-A.}}, \bauthor{\bsnm{Meinshausen},~\bfnm{Nicolai}\binits{N.}} \AND \bauthor{\bsnm{Shah},~\bfnm{Rajen~D}\binits{R.~D.}}
(\byear{2018}).
\btitle{The $xyz$ algorithm for fast interaction search in high-dimensional data}.
\bjournal{Journal of Machine Learning Research}
\bvolume{19}
\bpages{1343--1384}.
\end{barticle}
\endbibitem

\bibitem{tibshirani96}
\begin{barticle}[author]
\bauthor{\bsnm{Tibshirani},~\bfnm{Robert}\binits{R.}}
(\byear{1996}).
\btitle{Regression Shrinkage and Selection Via the Lasso}.
\bjournal{Journal of the Royal Statistical Society, Series B}
\bvolume{58}
\bpages{267--288}.
\end{barticle}
\endbibitem

\bibitem{valiant1984theory}
\begin{barticle}[author]
\bauthor{\bsnm{Valiant},~\bfnm{Leslie~G}\binits{L.~G.}}
(\byear{1984}).
\btitle{A theory of the learnable}.
\bjournal{Communications of the ACM}
\bvolume{27}
\bpages{1134--1142}.
\end{barticle}
\endbibitem

\bibitem{vanderweele2008empirical}
\begin{barticle}[author]
\bauthor{\bsnm{VanderWeele},~\bfnm{Tyler~J}\binits{T.~J.}} \AND \bauthor{\bsnm{Robins},~\bfnm{James~M}\binits{J.~M.}}
(\byear{2008}).
\btitle{Empirical and counterfactual conditions for sufficient cause interactions}.
\bjournal{Biometrika}
\bvolume{95}
\bpages{49--61}.
\end{barticle}
\endbibitem

\bibitem{vansteelandt2022assumption}
\begin{barticle}[author]
\bauthor{\bsnm{Vansteelandt},~\bfnm{Stijn}\binits{S.}} \AND \bauthor{\bsnm{Dukes},~\bfnm{Oliver}\binits{O.}}
(\byear{2022}).
\btitle{Assumption-lean inference for generalised linear model parameters}.
\bjournal{Journal of the Royal Statistical Society Series B: Statistical Methodology}
\bvolume{84}
\bpages{657--685}.
\end{barticle}
\endbibitem

\bibitem{Wahba1994}
\begin{barticle}[author]
\bauthor{\bsnm{Wahba},~\bfnm{Grace}\binits{G.}}, \bauthor{\bsnm{Wang},~\bfnm{Yuedong}\binits{Y.}}, \bauthor{\bsnm{Gu},~\bfnm{Chong}\binits{C.}}, \bauthor{\bsnm{Klein},~\bfnm{Ronald}\binits{R.}} \AND \bauthor{\bsnm{Klein},~\bfnm{Barbara}\binits{B.}}
(\byear{1995}).
\btitle{{Smoothing spline ANOVA for exponential families, with application to the Wisconsin Epidemiological Study of Diabetic Retinopathy: the 1994 Neyman Memorial Lecture}}.
\bjournal{The Annals of Statistics}
\bvolume{23}
\bpages{1865--1895}.
\bdoi{10.1214/aos/1034713638}
\end{barticle}
\endbibitem

\bibitem{wooldridge2010econometric}
\begin{bbook}[author]
\bauthor{\bsnm{Wooldridge},~\bfnm{Jeffrey~M}\binits{J.~M.}}
(\byear{2010}).
\btitle{Econometric analysis of cross section and panel data}.
\bpublisher{MIT press}.
\end{bbook}
\endbibitem

\bibitem{xiang2022pairwise}
\begin{barticle}[author]
\bauthor{\bsnm{Xiang},~\bfnm{Siqi}\binits{S.}}, \bauthor{\bsnm{Liu},~\bfnm{Siyao}\binits{S.}}, \bauthor{\bsnm{Perou},~\bfnm{Charles~M}\binits{C.~M.}}, \bauthor{\bsnm{Zhang},~\bfnm{Kai}\binits{K.}} \AND \bauthor{\bsnm{Marron},~\bfnm{JS}\binits{J.}}
(\byear{2022}).
\btitle{Pairwise Nonlinear Dependence Analysis of Genomic Data}.
\bjournal{arXiv preprint arXiv:2202.09880}.
\end{barticle}
\endbibitem

\bibitem{zhang2017spherical}
\begin{barticle}[author]
\bauthor{\bsnm{Zhang},~\bfnm{Kai}\binits{K.}}
(\byear{2017}).
\btitle{Spherical Cap Packing Asymptotics and Rank-Extreme Detection}.
\bjournal{IEEE Transactions on Information Theory}
\bvolume{63}
\bpages{4572-4584}.
\bdoi{10.1109/TIT.2017.2700202}
\end{barticle}
\endbibitem

\bibitem{zhang2019bet}
\begin{barticle}[author]
\bauthor{\bsnm{Zhang},~\bfnm{Kai}\binits{K.}}
(\byear{2019}).
\btitle{{BET} on Independence}.
\bjournal{Journal of the American Statistical Association}
\bvolume{114}
\bpages{1620-1637}.
\bdoi{10.1080/01621459.2018.1537921}
\end{barticle}
\endbibitem

\bibitem{zhang2021beauty}
\begin{barticle}[author]
\bauthor{\bsnm{Zhang},~\bfnm{Kai}\binits{K.}}, \bauthor{\bsnm{Zhang},~\bfnm{Wan}\binits{W.}}, \bauthor{\bsnm{Zhao},~\bfnm{Zhigen}\binits{Z.}} \AND \bauthor{\bsnm{Zhou},~\bfnm{Wen}\binits{W.}}
(\byear{2021}).
\btitle{{BEAUTY} powered {BEAST}}.
\bjournal{arXiv preprint arXiv:2103.00674}.
\end{barticle}
\endbibitem

\bibitem{zhao2021fast}
\begin{barticle}[author]
\bauthor{\bsnm{Zhao},~\bfnm{Zhigen}\binits{Z.}}, \bauthor{\bsnm{Zhang},~\bfnm{Wan}\binits{W.}}, \bauthor{\bsnm{Baiocchi},~\bfnm{Mike}\binits{M.}} \AND \bauthor{\bsnm{Zhang},~\bfnm{Kai}\binits{K.}}
(\byear{2021}).
\btitle{Fast, Flexible, and Powerful: Introducing a Scalable, Bitwise Framework for Non-parametric Testing for Dependence Structure}.
\bjournal{submitted}.
\end{barticle}
\endbibitem

\end{thebibliography}

\end{document}